\documentclass[11pt]{article}
\usepackage{amsfonts,mathrsfs,amssymb,amsthm,mathptm}
\usepackage{amsmath,amscd}
\usepackage{mathptm,pslatex}
\usepackage{multirow}
\usepackage{color}
\usepackage[dvips]{graphicx}
\usepackage[all]{xy}
\usepackage{graphicx}
\usepackage{fancyhdr}
\usepackage{url}
\usepackage{cite}
\usepackage{exscale}
\usepackage{relsize}
\usepackage{bbm}

\begin{document}
\newtheorem{Def}{Definition}[section]
\newtheorem{Bsp}[Def]{Example}
\newtheorem{Prop}[Def]{Proposition}
\newtheorem{Theo}[Def]{Theorem}
\newtheorem{Lem}[Def]{Lemma}
\newtheorem{Koro}[Def]{Corollary}
\newtheorem{Ques}[Def]{Question}
\theoremstyle{definition}
\newtheorem{Rem}[Def]{Remark}

\newcommand{\diag}{{\rm diag}}
\newcommand{\add}{{\rm add}}
\newcommand{\con}{{\rm con}}
\newcommand{\gd}{{\rm gl.dim}}
\newcommand{\dd}{{\rm dom.dim}}
\newcommand{\codd}{{\rm codom.dim}}
\def\ndd{\nu\mbox{-}\dd}
\newcommand{\tdim}{{\rm dim}}
\newcommand{\E}{{\rm E}}
\newcommand{\Mor}{{\rm Morph}}
\newcommand{\End}{{\rm End}}
\newcommand{\ind}{{\rm ind}}
\newcommand{\rsd}{{\rm res.dim}}
\newcommand{\rd} {{\rm rep.dim}}
\newcommand{\ol}{\overline}
\newcommand{\overpr}{$\hfill\square$}
\newcommand{\rad}{{\rm rad}}
\newcommand{\soc}{{\rm soc}}
\renewcommand{\top}{{\rm top}}
\newcommand{\stp}{{\mbox{\rm -stp}}}
\newcommand{\pd}{{\rm projdim}}
\newcommand{\id}{{\rm injdim}}
\newcommand{\del}{{\rm del}}
\newcommand{\fld}{{\rm flatdim}}
\newcommand{\fdd}{{\rm fdomdim}}
\newcommand{\Fac}{{\rm Fac}}
\newcommand{\Gen}{{\rm Gen}}
\newcommand{\fd} {{\rm findim}}
\newcommand{\Fd} {{\rm Findim}}
\newcommand{\Pf}[1]{{\mathscr P}^{<\infty}(#1)}
\newcommand{\DTr}{{\rm DTr}}
\newcommand{\cpx}[1]{#1^{\bullet}}
\newcommand{\D}[1]{{\mathscr D}(#1)}
\newcommand{\Dz}[1]{{\mathscr D}^+(#1)}
\newcommand{\Df}[1]{{\mathscr D}^-(#1)}
\newcommand{\Db}[1]{{\mathscr D}^b(#1)}
\newcommand{\C}[1]{{\mathscr C}(#1)}
\newcommand{\Cz}[1]{{\mathscr C}^+(#1)}
\newcommand{\Cf}[1]{{\mathscr C}^-(#1)}
\newcommand{\Cb}[1]{{\mathscr C}^b(#1)}
\newcommand{\Dc}[1]{{\mathscr D}^c(#1)}
\newcommand{\K}[1]{{\mathscr K}(#1)}
\newcommand{\Kz}[1]{{\mathscr K}^+(#1)}
\newcommand{\Kf}[1]{{\mathscr  K}^-(#1)}
\newcommand{\Kb}[1]{{\mathscr K}^b(#1)}
\newcommand{\modcat}{\ensuremath{\mbox{{\rm -mod}}}}
\newcommand{\Modcat}{\ensuremath{\mbox{{\rm -Mod}}}}
\newcommand{\stmodcat}[1]{#1\mbox{{\rm -{\underline{mod}}}}}
\newcommand{\pmodcat}[1]{#1\mbox{{\rm -proj}}}
\newcommand{\imodcat}[1]{#1\mbox{{\rm -inj}}}
\newcommand{\Pmodcat}[1]{#1\mbox{{\rm -Proj}}}
\newcommand{\Imodcat}[1]{#1\mbox{{\rm -Inj}}}
\newcommand{\PI}[1]{#1\mbox{{\rm -prinj}\,}}
\newcommand{\opp}{^{\rm op}}
\newcommand{\otimesL}{\otimes^{\rm\mathbb L}}
\newcommand{\rHom}{{\rm\mathbb R}{\rm Hom}\,}
\newcommand{\projdim}{\pd}
\newcommand{\Hom}{{\rm Hom}}
\newcommand{\Coker}{{\rm Coker}}
\newcommand{ \Ker  }{{\rm Ker}}
\newcommand{ \Cone }{{\rm Con}}
\newcommand{ \Img  }{{\rm Im}}
\newcommand{\Ext}{{\rm Ext}}
\newcommand{\StHom}{{\rm \underline{Hom}}}
	
\newcommand{\gm}{{\rm _{\Gamma_M}}}
\newcommand{\gmr}{{\rm _{\Gamma_M^R}}}
	
\def\vez{\varepsilon}\def\bz{\bigoplus}  \def\sz {\oplus}
\def\epa{\xrightarrow} \def\inja{\hookrightarrow}
	
\newcommand{\lra}{\longrightarrow}
\newcommand{\llra}{\longleftarrow}
\newcommand{\lraf}[1]{\stackrel{#1}{\lra}}
\newcommand{\llaf}[1]{\stackrel{#1}{\llra}}
\newcommand{\ra}{\rightarrow}
\newcommand{\dk}{{\rm dim_{_{k}}}}
\newcommand{\colim}{{\rm colim\, }}
\newcommand{\limt}{{\rm lim\, }}
\newcommand{\Add}{{\rm Add }}
\newcommand{\Tor}{{\rm Tor}}
\newcommand{\Cogen}{{\rm Cogen}}
\newcommand{\Tria}{{\rm Tria}}
\newcommand{\tria}{{\rm tria}}

{\Large \bf
\begin{center}
New invariants of stable equivalences and Auslander--Reiten conjecture
\end{center}}
	
\medskip
\centerline{\textbf{Changchang Xi and Jinbi Zhang$^*$}}

\medskip
\centerline{In memory of Professor Idun Reiten (1942--2025)}
	
\renewcommand{\thefootnote}{\alph{footnote}}
\setcounter{footnote}{-1} \footnote{$^*$Corresponding author's Email: zhangjb@cnu.edu.cn}
\renewcommand{\thefootnote}{\alph{footnote}}
\setcounter{footnote}{-1} 
\renewcommand{\thefootnote}{\alph{footnote}}
\setcounter{footnote}{-1} 

\begin{abstract}
We show that stable equivalences between Artin algebras without nodes preserve homological data that provide upper bounds for finitistic dimension, and that stable equivalences between  Artin algebras with positive $\nu$-dominant dimensions induce stable equivalences of their Frobenius parts. As an application of our new methods developed, we verify the Auslander--Reiten conjecture on stable equivalences for two rather different classes of algebras: principal centralizer matrix algebras over arbitrary fields and Frobenius-finite algebras over algebraically closed fields.
\end{abstract}

{\footnotesize\tableofcontents\label{contents}}

\section{Introduction}
In the representation theory of finite-dimensional algebras and finite groups, there are several unsolved major conjectures that are of great interest. Two of them are the famous Auslander--Reiten conjecture on stable equivalences (ARC) (see \cite[Conjecture (5), p.409]{Aus1997}, or \cite[Conjecture 2.5]{rouquier}), and the finitistic dimension conjecture (FDC) (see \cite{bass}, or \cite[Conjecture (11), p.410]{Aus1997}).

(ARC) Stably equivalent Artin algebras have the same number of non-isomorphic, non-projective simple modules.

(FDC) Every Artin algebra has finite finitistic dimension.

Both conjectures are known in few cases only. For instance,  Auslander and Reiten proved that (ARC) holds for a stable equivalence between an Artin algebra and a hereditary algebra \cite{AR1975}, Mart\'inez-Villa showed that (ARC) holds for stable equivalences between representation-finite algebras over an algebraically closed field \cite{MV1985}, and reduced the validity of (ARC) to self-injective algebras without nodes (see \cite{MV1990}). Recently, the conjecture was proved for stable equivalences of Morita type between Frobenius-finite algebras over algebraically closed fields and without semisimple summands (see \cite{hx18}), and for stable equivalences between special biserial algebras (see \cite{AnZv2017},\cite{Pogorzaly1994}). As to (FDC), there are several approaches. For example, it was verified for algebras of radical cube-zero \cite{green-zim} and monomial algebras \cite{gkirkman}. Further, Igusa and Todorov introduced two homological data, called $\phi$- and $\psi$-dimensions, to give upper bounds of the finitistic dimensions \cite{it05}. Another approach to (FDC) was given for algebras over perfect fields by relative homological algebra \cite{xx}. Recently, G\'elinas has introduced a new homological datum, called the delooping level of an algebra in \cite{Gelinas22}, which provides upper bounds of the finitistic dimensions of opposite algebras (see \cite[Proposition 1.3]{Gelinas22}).

In this article, we provide new stable equivalences and new invariants of stable equivalences. Further, we use these to prove the validity of the conjectures or to transfer validity from algebras to stably equivalent algebras.

\begin{Theo}\label{arc+del}
Stably equivalent Artin algebras without nodes preserve delooping level as well as Igusa and Todorov's $\phi$-dimension and $\psi$-dimension.
\end{Theo}

As a consequence, if one of the algebras has finite $\phi$- or $\psi$-dimensions, then both algebras have finite finitistic dimensions. Moreover, if the opposite algebra of one of the algebras has finite delooping level, then both algebras have finite
finitistic dimension. For the definitions of $\phi$- and $\psi$-dimensions as well as delooping levels, we refer to Subsection \ref{sect2.2} below.

Now, we apply Theorem \ref{arc+del} to \emph{Morita algebras} that are just the endomorphism algebras of modules of the form $A\oplus M$ with $A$ a self-injective algebra and $M$ a finitely generated module over $A$ (see \cite{KY13}). This class of algebras includes many algebras arising in the representation theory or algebraic Lie theory such as Temperley--Lieb algebras, Schur algebras and $q$-Schur algebras. Our application can be stated in a general form which involves almost $\nu$-stable derived equivalences. The latter induces stable equivalences.

\begin{Theo}
\label{a-n-s-d-e}
$(1)$ Almost $\nu$-stable derived equivalences between finite-dimensional algebras over a field preserve delooping level as well as Igusa and Todorov's $\phi$-dimension and $\psi$-dimension.

$(2)$ If $A$ is a self-injective algebra over a field and $X$ is a finitely generated $A$-module, then
all of the endomorphism algebras of the modules $A\oplus X, A\oplus \Omega_A(X)$ and $A\oplus D{\rm Tr}(X)$ have the same delooping level as well as the same Igusa and Todorov's $\phi$-dimension and $\psi$-dimension, where $\Omega_A$ and $D{\rm Tr}$ are the syzygy operator and Auslander--Reiten translation of $A$, respectively.
\end{Theo}

As is known, the Auslander--Reiten conjecture can be reduced to stable equivalences between self-injective algebras \cite{MV1990}. Thus it is of great importance to know whether a stable equivalence induces the one between their Frobenius parts. On the other hand, it is a fundamental, but difficult problem to get examples of stable equivalences between algebras. In the literature, only a few classes of stable equivalences have been constructed (for example, see \cite{AR1975, Reiten78, hx10, lx05}). The next result provides a rather large class of new stable equivalences.

\begin{Theo}\label{main-thm}
Stable equivalences between Artin algebras of positive $\nu$-dominant dimensions induce the ones between their Frobenius parts.
\end{Theo}

Compared with a result \cite[Theorem 17]{MV1990}, Theorem \ref{main-thm} does not assume Artin algebras considered to have neither nodes nor semisimple direct summands.

\medskip
Finally, as an application, we verify the Auslander--Reiten conjecture on stable equivalences for two new classes of algebras, namely Frobenius-finite algebras and principal centralizer matrix algebras. Recall that an Artin algebra is said to be \emph{Frobenius-finite} if its Frobenius part is representation-finite. For the definition of Frobenius parts of algebras, we refer to Subsection \ref{sect2.3} below or \cite{hx18, MV1985}. A \emph{principal centralizer matrix algebra} is the centralizer of an $n\times n$ matrix in the full $n\times n$ matrix algebra over a field.
Though these two large classes of algebras are rather different in nature, our approach to the Auslander--Reiten conjecture works for both of them.

\begin{Theo} \label{cma} The Auslander--Reiten conjecture on stable equivalences holds true for the two classes of finite-dimensional algebras:

$(1)$ Principal centralizer matrix algebras over a field.

$(2)$ Frobenius-finite algebras over an algebraically closed field.
\end{Theo}

If we assume in Theorem \ref{cma}(2) that the stable equivalences are of Morita type and both algebras have no semisimple direct summands, then Theorem \ref{cma}(2) follows from \cite[Theorem 1.1]{hx18}. Since the class of Frobenius-finite algebras properly contains the one of representation-finite algebras, Theorem \ref{cma}(2) generalizes also a result in \cite[Theorem 3.4]{MV1985} which states that the Auslander--Reiten conjecture on stable equivalences holds true for representation-finite algebras over an algebraically closed field.

\medskip
This paper is structured as follows. In Section \ref{prel}, we recall definition of stable equivalences, finitistic dimensions, almost $\nu$-stable equivalences and homological data. In Section \ref{sect3}, we give necessary information and prove some general results on stable equivalences for later proofs. In Section \ref{proof-main}, we prove Theorems \ref{arc+del}-\ref{main-thm}. We point out that tilting does not preserve delooping levels. In Section \ref{sect5}, we apply our methods to prove Theorem \ref{cma}. We conjecture that finiteness of delooping levels of Artin algebras is invariant under derived equivalences.

\medskip
{\bf Acknowledgements.} The research work was partially supported by the National Natural Science Foundation of China (Grants 12031014 and 12401038). The authors thank Dr. Xiaogang Li for some discussions on centralizer matrix algebras.

Both authors are deeply grateful to the editor and the anonymous reviewer for carefully reading the manuscript and proposing a lot of constructive suggestions and helpful comments for improving the presentation of results in this article.

\section{Preliminaries}\label{prel}
In this section, we fix notations and recall several homological data as well as basic results of stable equivalences.

\subsection{Auslander--Reiten conjecture and finitistic dimension conjecture}
Let $A$ be an Artin algebra over a commutative Artin ring $k$. By $A\modcat$ we denote the category of all finitely generated left $A$-modules. Let $A^{\opp}$ be the opposite algebra of $A$, we understand a right $A$-module as a left $A^{\opp}$-module. We denote by $D$ the usual duality of Artin algebra from $A\modcat$ to $A^{\opp}\modcat$. For $M\in A\modcat$, let $\Omega_A(M)$ be the syzygy of $_AM$; Tr$(M)$ the transpose of $M$, which is an $A^{\opp}$-module; and $\add(M)$  the full additive subcategory of $A\modcat$ consisting of all direct summands of finite sums of copies of $M$. We write $\End_A(M)$ for the endomorphism algebra of $_AM$.

Let $A\modcat_{\mathscr{P}}$ (respectively, $A\modcat_{\mathscr{I}}$)  be the full subcategory of $A\modcat$ consisting of those modules that do not have nonzero projective (respectively, injective) direct summands. Let $\mathscr{P}(A)_{\mathscr{I}}$ (respectively, $\mathscr{I}(A)_{\mathscr{P}}$) denote the set of all isomorphism classes of indecomposable projective (respectively, injective) $A$-modules without any nonzero injective (respectively, projective) direct summands.

The stable category $A\stmodcat$ of $A$ has the same objects as $A\modcat$, its morphism set $\underline{\Hom}_A(X,Y)$ of two modules $X$ and $Y$ is the quotient $k$-module of $\Hom_A(X,Y)$ modulo all homomorphisms that factorize through a projective $A$-module. For $f\in\Hom_A(X,Y)$, we write $\underline{f}$ for the image of $f$ in $\StHom_A(X,Y)$. Note that $X\simeq Y$ in $A\stmodcat$ if and only if there are two projective modules $P, Q\in A\modcat$ such that $X\oplus P\simeq Y\oplus Q$ as $A$-modules. In this case, $\End_A(A\oplus X)$ and $\End_A(A\oplus Y)$ are Morita equivalent.

Let $(\stmodcat{A})\modcat$ be the category of all finitely presented functors from $(\stmodcat{A})^{\rm op}$ to the category $\mathscr{A}$ of all abelian groups. Recall that a functor $H:(\stmodcat{A})^{\rm op} \ra \mathscr{A}$ is said to be \emph{finitely presented} if there exists an exact sequence of functors $\StHom_A(-,X)\ra \StHom_A(-,Y)\ra H\ra 0$ with $X$ and $Y$ in $\stmodcat{A}$. It is known from \cite{ar73} that $(\stmodcat{A})\modcat$ is an abelian category, its projective objects are precisely the functors $\StHom_A(-,X)$ for $ X\in A\modcat_{\mathscr{P}}$, and its injective objects are precisely the functors $\Ext^1_A(-,X)$ for $ X\in A\modcat_{\mathscr{P}}$.

Artin algebras $A$ and $B$ over a commutative Artin ring $k$ are said to be \emph{stably equivalent} if the two stable categories $A\stmodcat$ and $B\stmodcat$ are equivalent as $k$-categories.

The \emph{Auslander--Reiten conjecture} on stable equivalences says that the number of non-isomorphic, non-projective simple modules is an invariant of stable equivalences. This is also called Auslander--Alperin conjecture in the representation theory of groups (see \cite{rouquier}).

By eliminating nodes and considering special projective-injective modules, Mart\'inez-Villa reduced this conjecture to the one just for self-injective Artin algebras \cite{MV1990}.
Following \cite{MV1980}, a non-projective, non-injective simple $A$-module $S$ is called a \emph{node} of $A$ if the middle term $P$ of the almost split sequence starting at $S$, $0\ra S\ra P\ra {\rm Tr\,}D(S)\ra 0$, is projective.

For an Artin algebra $A$, the \emph{finitistic dimension} of $A$, denoted $\fd(A)$, is the supremum
of the projective dimensions of modules $X\in A\modcat$ that have finite projective dimension. The \emph{finitistic dimension conjecture} says that every Artin algebra has finite finitistic dimension (see \cite{bass}).

Considerable efforts notwithstanding, the above two conjectures are still open to date. In the next subsection, we will mention two approaches to the finitistic dimension conjecture by Igusa--Todorov, and G\'elinas.

\subsection{Homological data: delooping levels and Igusa-Todorov's dimensions\label{sect2.2}}
In this subsection, we recall a few homological data: delooping levels, $\phi$- and $\psi$-dimensions.

To understand the finitistic dimensions of algebras, Igusa and Todorov introduced the $\phi$- and $\psi$-dimensions for Artin algebras in \cite{it05}.

Let $K(A)$ be the Grothendieck group of $A$, that is, the quotient of the free abelian group generated by the isomorphism classes $[X]$ with $X\in A\modcat$, modulo the relations:

(1) $[Z]=[X]+[Y]$ if $_AZ\simeq {}_AX\oplus {}_AY$;

(2) $[P]=0$ if $_AP$ is projective.

\noindent Then $K(A)$ is the free abelian group generated by the isomorphism classes of non-projective indecomposable $A$-modules $X\in A\modcat$. Now, we recall two functions $\phi$ and $\psi$ from $A\modcat$ to $\mathbb{N}$, the set of natural numbers, defined in \cite{it05}.

The syzygy functor $\Omega:\stmodcat{A} \ra A\stmodcat$ on the stable module category $A\stmodcat$ induces a group homomorphism $\Omega:K(A)\ra K(A)$ of abelian groups, given by $\Omega([X]):=[\Omega(X)]$. For $X\in A\modcat$, let $\langle X\rangle$ be the $\mathbb{Z}$-submodule of $K(A)$ generated by the isomorphism classes of non-projective, indecomposable direct summands of $X$. Since the
rank of the image $\Omega\langle X\rangle$ of $\langle X\rangle$ under $\Omega$ does not exceed the finite rank of $\langle X\rangle$, it follows from Fitting's Lemma that there exists a smallest nonnegative integer $\phi(X)$ such that $\Omega:\Omega^n\langle X\rangle\ra \Omega^{n+1}\langle X\rangle$ is an isomorphism for all $n\ge \phi(X)$. Furthermore, let
$$\psi(X):=\phi(X)+\sup\{\pd(Y)\mid Y \mbox{ is a direct summand of }\Omega^{\phi(X)}(X), \; \pd(Y)<\infty\}.$$
Then the $\phi$-dimension and $\psi$-dimension of $A$ are defined by
$$\phi\dim(A):=\sup\{\phi(X)\mid X\in A\modcat\}\quad \mbox{ and } \quad \psi\dim(A):=\sup\{\psi(X)\mid X\in A\modcat\}.$$
According to \cite[Lemma 0.3]{it05}, $\psi(X)=\phi(X)=\pd(X)$ if $\pd(X)<\infty$.  Thus $$\fd(A)\le\phi\dim(A)\le \psi\dim(A)\le\gd(A),$$where $\gd(A)$ means the global dimension of $A$.

Recently, G\'elinas has introduced a new homological datum: the delooping level of an algebra in \cite{Gelinas22}. The significance of delooping levels is that the finitistic dimensions of algebras can be bounded by the delooping levels of the opposite algebras (see \cite[Proposition 1.3]{Gelinas22}).

By definition, the \emph{delooping level} of an $A$-module $X\in A\modcat$, denoted $\del(X)$, is the smallest number $d\ge 0$ such that the $d$-th syzygy $\Omega^d(X)$ of $X$ is a direct summand of a module of the form $P\oplus \Omega^{d+1}(M)$ for an $A$-module $M$ and a projective $A$-module $P$, where $\Omega$ is the syzygy (or loop) operator of $A$. If such a number $d$ does not exist, one defines $\del(X)=\infty$. The \emph{delooping level} of $A$, denoted $\del(A)$, is the maximum of the delooping levels of all non-isomorphic simple $A$-modules. Clearly, $\del(X\oplus Y)=\max\{\del(X),\del(Y)\}$ for $X,Y\in A\modcat$,
and $\del(A)=\del(\top(_AA))$, where $\top(_AX)$ denotes the top of an $A$-module $_AX$. It was shown in \cite{Gelinas22} that $\fd(A)\le \del(A^{\opp})$, where $A^{\opp}$ stands for the opposite algebra of $A$.

As is known, each of these homological data may be infinite. This implies that the finitistic dimension conjecture is still open to date.

\subsection{Frobenius parts and $\nu$-dominant dimensions\label{sect2.3}}

Let $A$ be an Artin algebra. We denote by $A$-prinj the full subcategory of $A\modcat$ consisting of those $A$-modules that are both projective and injective.

A projective $A$-module $P$ is said to be \emph{$\nu$-stably projective} \cite{hx18} if $\nu_A^i(P)$ is projective for all $i> 0$, where $\nu_A$ is the Nakayama functor of $A$. By $A$-stp we denote the full subcategory of $A\modcat$ consisting of all $\nu$-stably projective $A$-modules. Clearly, $A\mbox{-stp} \subseteq A\mbox{-prinj}$.

If $X\in A\modcat$ such that $\add(X)=A\stp$, then the endomorphism algebra $\End_A(X)$ of $X$ is called a \emph{Frobenius part} of $A$. It is a self-injective algebra (see \cite{MV1985}, or \cite[Lemma 2.7]{hx18}) and is defined uniquely up to Morita equivalence. Note that Frobenius parts of Artin algebras were first given by Mart\'inez-Villa in different but equivalent terms in \cite{MV1985}. Following \cite{hx18}, an Artin algebra is said to be \emph{Frobenius-finite} if its Frobenius part is representation-finite.

For an $A$-module $M\in A\modcat$, we consider its minimal injective resolution
$$0\lra {}_AM\lra I_0\lra I_1\lra I_2\lra\cdots.$$
Let $I$ be an injective $A$-module and $0\le n\le \infty$. If $n$ is maximal such that all modules $I_j$ are in $\add(I)$ for $j < n$, then $n$ is called the \emph{$I$-dominant dimension} of $M$, denoted by $I$-$\dd(M)$. Dually, we consider its minimal projective resolution
$$\cdots\lra P_2\lra P_1\lra P_0\lra {}_AM \lra 0.$$
Let $P$ be a projective $A$-module and $0\le m\le \infty$. If $m$ is maximal such that all modules $P_j$ are in $\add(P)$ for $j < m$, then $m$ is called the \emph{$P$-codominant dimension} of $M$, denoted by $P$-$\codd(M)$. Clearly, the codominant dimension of $M$ is equal to the dominant dimension of $A^{\opp}$-module  $D(M)$.
Now, if $\add(I)=\add(P) = \PI{A}$, then we define the \emph{dominant dimension} of $M$ to be $I$-$\dd(M)$, denoted by $\dd(M)$; the \emph{codominant dimension} of $M$ to be  $P$-$\codd(M)$, denoted by $\codd(M)$, and the \emph{dominant dimension} of the algebra $A$ to be $\dd(_AA)$, denoted by $\dd(A)$. Note that $\dd(A) =\dd(A^{\opp})$ (see \cite[Theorem 4]{M68} or \cite{Hos89}). It is clear that $\dd(A) = \min\{\dd(P)\mid P\in \add(_AA)\}$.
If $\add(I)=A\stp$, then $I$-$\dd(M)$ is called the $\nu$-\emph{dominant dimension} of $M$, denoted by $\ndd(M)$. The $\nu$-\emph{dominant dimension of the algebra $A$} is defined to be $\ndd(_AA)$.

\begin{Lem}\label{pi=stp}
Let $A$ be an Artin algebra with $\ndd(A)\ge 1$. Then

$(1)$ $A\stp=\PI{A}$ and $\ndd(A)=\dd(A)$.

$(2)$ The projective cover of a simple module $_AS$ is injective if and only if the injective envelope of $S$ is projective.

$(3)$ If the projective cover of a simple module $_AS$ is not injective, then $S$ itself is neither projective nor injective.
\end{Lem}
{\it Proof.} (1) and (2) are trivial. We prove (3). Let $P$ and $I$ be the projective cover and injective envelope of $S$, respectively. By assumption, $P\not\in \PI{A}$. By (2), $I\not\in \PI{A}$. If $S$ is injective, then it follows from $\codd(D(A_A))=\dd(A)=\ndd(A)\ge 1$ that $\codd(S)\ge 1$. This implies $P\in \PI{A}$, a contradiction. Thus $S$ is not an injective module. Suppose that $S$ is projective. It follows from $\dd(A)\ge 1$ that $\dd(S)\ge 1$ and $I\in \PI{A}$, again a contradiction. Thus $S$ is not a projective module.
$\square$

\section{Basics of stable equivalences\label{sect3}}
In this section, we recall basic facts on stable equivalences and develop some lemmas for later proofs.

\subsection{General facts on stable equivalences\label{sect3.1}}
In this subsection we collect some known facts and prove a new lemma on general stable equivalences.

\medskip
Assume that $F: \stmodcat{A} \ra \stmodcat{B}$ is an equivalence of $k$-categories with a quasi-inverse functor $G: \stmodcat{B} \ra \stmodcat{A}$. So $F$ and $G$ are additive functors and induce two equivalences $\alpha$ and $\beta$ of abelian categories (see \cite[Section 8]{ar74})
$$\alpha:(\stmodcat{A})\modcat\lraf{\simeq} (\stmodcat{B})\modcat\quad \mbox{and}\quad \beta:(\stmodcat{B})\modcat\lraf{\simeq} (\stmodcat{A})\modcat,$$
and two one-to-one correspondences
$$F: A\modcat_{\mathscr{P}}\longleftrightarrow B\modcat_{\mathscr{P}}:G \quad \mbox{and}\quad F': A\modcat_{\mathscr{I}}\longleftrightarrow  B\modcat_{\mathscr{I}}:G'$$
such that $$\alpha(\StHom_A(-,X))\simeq \StHom_B(-,F(X))\quad \mbox{and} \quad \alpha(\Ext_A^1(-,Y))\simeq \Ext_B^1(-,F'(Y)),$$
$$\beta(\StHom_B(-,U))\simeq \StHom_A(-,G(U))\quad \mbox{and} \quad \beta(\Ext_B^1(-,V))\simeq \Ext_A^1(-,G'(V))$$
for $X\in A\modcat_{\mathscr{P}},\, Y\in A\modcat_{\mathscr{I}}$, $U\in B\modcat_{\mathscr{P}}$ and $V\in B\modcat_{\mathscr{I}}$. For convenience, we set $F(P)=0$ for a projective module $P$, and $F'(I)=0$ for an injective module $I$.

The following lemma is useful for later discussions.
\begin{Lem}\label{ext-iso}{\rm (\cite[Section 7, p.347]{ar74})}
If $X,Y\in A\modcat_{\mathscr{I}}$, then $X\simeq Y$ in $A\modcat$ if and only if $\Ext_A^1(-,X)\simeq \Ext^1_A(-,Y)$ in $(\stmodcat{A})\modcat$.
\end{Lem}

A node $S$ of $A$ is called an $F$-\emph{exceptional node} if $F(S)\not\simeq F'(S)$. By $\mathfrak{n}_{F}(A)$ we denote the set of isomorphism classes of $F$-exceptional nodes of $A$. By \cite[Lemma 3.4]{ar78}, if an indecomposable $A$-module $X$ is neither a node, nor projective and nor injective, then $F(X)\simeq F'(X)$. Thus $\mathfrak{n}_{F}(A)$ coincides with the set of isomorphism classes
of non-projective, non-injective, indecomposable $A$-modules $X$ such that $F(X)\not\simeq F'(X)$.

In the following, let $$\bigtriangleup_A
:=\mathfrak{n}_{F}(A)
\dot{\cup}\mathscr{P}(A)_{\mathscr{I}}
\quad\mbox{and}\quad \bigtriangledown_A
:=\mathfrak{n}_{F}(A)
\dot{\cup}\mathscr{I}(A)_{\mathscr{P}},$$ where $\dot{\cup}$ stands for the disjoint union of sets; $\mathscr{P}(A)_{\mathscr{I}}$ (respectively, $\mathscr{I}(A)_{\mathscr{P}}$) stands for the set of representatives of isomorphism classes of indecomposable projective (respectively, injective) $A$-modules without any nonzero injective (respectively, projective) summands. By $\bigtriangleup_A^c$ we mean the class of non-injective, indecomposable $A$-modules which do not belong to $\bigtriangleup_A$. Thus
each module $M\in A\modcat_{\mathscr{I}}$ admits a unique decomposition (up to isomorphism)
$$M\simeq M_1\oplus M_2$$
with $M_1\in \add(\bigtriangleup_A)$ and $M_2\in \add(\bigtriangleup_A^c)$.

\begin{Lem}\label{one-to-one}
{\rm (\cite[Lemma 4.10(1)]{c21})}
The functor $F$ induces the bijections
$$F:\bigtriangledown_A\longleftrightarrow \bigtriangledown_B:G
,\quad
F':\bigtriangleup_A\longleftrightarrow \bigtriangleup_B:G'
\quad\mbox{and}\quad
F':\bigtriangleup_A^c\longleftrightarrow \bigtriangleup_B^c:G'.$$
\end{Lem}

An exact sequence $0\ra X\stackrel{f}{\ra} Y\stackrel{g}{\ra} Z\ra 0$ in $A\modcat$ is called \emph{minimal} if it does not have a split exact sequence as its direct summand, that is, there do not exist isomorphisms $u$, $v$, $w$ such that the diagram
$$\xymatrix{
0\ar[r]
&X\ar[r]^{f}\ar[d]^{u}
&Y\ar[r]^g\ar[d]^{v}
&Z\ar[r]\ar[d]^{w}
&0\\
0\ar[r]&X_1\oplus
X_2\ar[r]^{\text{
$\left(
\begin{smallmatrix}
f_1&0\\
0&f_2
\end{smallmatrix}\right)
$}}
&Y_1\oplus Y_2\ar[r]^{\text{
$\left(
\begin{smallmatrix}
g_1&0\\
0&g_2
\end{smallmatrix}\right)
$}}&Z_1\oplus Z_2\ar[r]&0
}$$
is commutative and exact in $A\modcat$, where $Y_2\ne 0$ and the sequence $0\ra X_2\lraf{f_2} Y_2\lraf{g_2}Z_2\ra 0$ splits. By definition, a minimal exact sequence does not split.

\begin{Lem}\label{lem-min}{\rm (\cite[Theorem 7.5]{ar74} or \cite[Proposition 2.1]{ar78})}
Let $H\in (\stmodcat{A})\modcat$ and $0\ra X\ra Y\ra Z\ra 0$
be a minimal exact sequence in $A\modcat$ such that the induced sequence $$0\lra \StHom_A(-,X)\lra\StHom_A(-,Y)\lra \StHom_A(-,Z)\lra H\lra 0$$
of functors is exact.  Then the following hold.

$(1)$ The induced exact sequence of functors
$$\StHom_A(-,Y)\lra \StHom_A(-,Z)\lra H\lra 0$$
is a minimal projective presentation of $H$ in $(\stmodcat{A})\modcat$.

$(2)$ The induced exact sequence of functors
$$0\lra H\lra \Ext_A^1(-,X)\lra \Ext_A^1(-,Y)$$
is a minimal injective copresentation of $H$ in $(\stmodcat{A})\modcat$.
\end{Lem}

The following lemma is from \cite[Lemma 1.6]{MV1990}, while its proof is referred to \cite{ar78}.
\begin{Lem}\label{minimal}
If $H\in (\stmodcat{A})\modcat$ has a minimal projective presentation
$$\xymatrix{\StHom_A(-,Y)\;\;\ar[r]^-{\StHom_A(-,g)}\ar[r] & \;\;\StHom_A(-,Z) \longrightarrow
 H \longrightarrow 0}$$
with $Y,Z\in A\modcat_{\mathscr{P}}$, then there is a minimal exact sequence
$$0\lra X\lra Y\oplus P\lraf{
g'} Z\lra 0$$
in $A\modcat$, where $\underline{g'}=\underline{g}$ in $\stmodcat{A}$ and $P$ is a projective $A$-module.
\end{Lem}

The following is a generalization of \cite[Theorem 1.7]{MV1990} and shows that the functor $F$ possesses certain ``exactness" property.
\begin{Lem}\label{exact-seq}
Let $0\ra X\oplus X'\ra Y\oplus \bar{Y}\oplus I\oplus P\oplus P'\stackrel{g}{\ra}Z\ra 0$
be a minimal exact sequence in $A\modcat$ with $X,Y\in \add(\bigtriangleup_A^c)$, $X'\in \add(\bigtriangleup_A)$, $\bar{Y}\in \add(\mathfrak{n}_{F}(A))$, $I\in\add(\mathscr{I}(A)_{\mathscr{P}})$, $P\in\add(\mathscr{P}(A)_{\mathscr{I}})$, $P'\in \PI{A}$ and $Z\in A\modcat_{\mathscr{P}}$.
Then there exists a minimal exact sequence
$$0\lra F(X)\oplus F'(X')\lra F(Y\oplus \bar{Y} \oplus I)\oplus Q\oplus Q'\lraf{g'}F(Z)\lra 0$$
in $B\modcat$, where $Q$ lies in $\add(\mathscr{P}(B)_{\mathscr{I}})$ and $Q'$ belongs to $\PI{B}$ such that $F(\bar{Y}\oplus I)\oplus Q \simeq F'(\bar{Y}\oplus P)\oplus J$ for some $J\in \add(\mathscr{I}(B)_{\mathscr{P}})$ and $\underline{g}'=F(\underline{g})$ in $\stmodcat{B}$.
\end{Lem}
{\it Proof.} We provide a proof by using some idea in \cite[Lemma 4.13]{c21}. Consider the finitely presented functor $H$:
$$ \StHom_A(-,Y\oplus \bar{Y}\oplus I\oplus P\oplus P')\lra \StHom_A(-,Z)\lra H\lra 0$$induced from the given minimal exact sequence
$$ \qquad 0\lra X\oplus X'\lra Y\oplus \bar{Y}\oplus I\oplus P\oplus P'\lraf{g}Z\lra 0$$ in $A\modcat$
with $I\in\add(\mathscr{I}(A)_{\mathscr{P}})$, $P\in\add(\mathscr{P}(A)_{\mathscr{I}})$ and $P'\in \PI{A}$. It follows from
Lemma \ref{lem-min}
that the sequence of functors
$$\xymatrix{\StHom_A(-,Y\oplus \bar{Y}\oplus I)\quad \ar[r]^-{{\tiny \StHom_A(-,g)}} &\quad\StHom_A(-,Z)\longrightarrow H\longrightarrow 0}$$
is a minimal projective presentation of $H$ in $(\stmodcat{A})\modcat$ and
that the sequence of functors
$$ 0\lra H\lra \Ext_A^1(-,X\oplus X')\lra \Ext_A^1(-,Y\oplus \bar{Y}\oplus P)$$ is a minimal injective copresentation of $H$ in $(\stmodcat{A})\modcat$.
Applying the equivalence functor $\alpha$ to the above two sequences of functors, we see that
the sequence
$$(\star)\qquad\StHom_B(-,F(Y\oplus \bar{Y}\oplus I)) \lra\StHom_B(-,F(Z))\lra \alpha(H) \lra 0$$
is a minimal projective presentation of $\alpha(H)$ in $(\stmodcat{B})\modcat$ and
the sequence
$$(*)\qquad 0\lra \alpha(H)\lra \Ext_B^1(-,F'(X)\oplus F'(X'))\lra \Ext_B^1(-,F'(Y \oplus \bar{Y} \oplus P))$$ is a minimal injective copresentation of $\alpha(H)$ in $(\stmodcat{B})\modcat$.

It follows from $(\star)$ and Lemma \ref{minimal} that there is a minimal exact sequence
$$(\diamondsuit)\qquad 0\lra W\lra F(Y\oplus \bar{Y} \oplus I)\oplus Q\oplus Q'\lraf{g'}F(Z)\lra 0$$
in $B\modcat$ with $Q\in\add(\mathscr{P}(B)_{\mathscr{I}})$, $Q'\in \PI{B}$ and $\underline{g}'=F(\underline{g})$ in $\stmodcat{B}$. The minimality of this sequence implies $W\in B\modcat_{\mathscr{I}}$. Note that $(\star)$ is induced from $(\diamondsuit)$.
Now, by Lemma \ref{lem-min}(2) and $(\diamondsuit)$, the exact sequence
$$ (\dag)\qquad 0\lra \alpha(H)\lra \Ext_B^1(-,W)\lra \Ext_B^1(-,F(Y\oplus \bar{Y} \oplus I)\oplus Q)$$of functors is a minimal injective copresentation of $\alpha(H)$ in $(\stmodcat{B})\modcat$. Thus both $(\dag)$ and $(*)$ are minimal injective copresentations of $\alpha(H)$. This implies that
$$(**)\qquad\Ext_B^1(-,F'(X)\oplus F'(X'))\simeq \Ext_B^1(-,W)\;
\mbox{ and}$$
$$\qquad\; \qquad (\ddag)\qquad\Ext_B^1(-,F'(Y \oplus \bar{Y} \oplus P))\simeq \Ext_B^1(-,F(Y\oplus \bar{Y} \oplus I)\oplus Q)$$
in $(\stmodcat{B})\modcat$. Since $X$ lies in $\add(\bigtriangleup_A^c)$ and $X'$ lies in $\add(\bigtriangleup_A)$, we know from Lemma \ref{one-to-one} that $F'(X)\in \add(\bigtriangleup_B^c)$ and $F'(X')\in \add(\bigtriangleup_B)$. In particular, $F'(X)\oplus F'(X')\in B\modcat_{\mathscr{I}}$. Thus $F'(X)\oplus F'(X')\simeq W$ as $B$-modules by Lemma \ref{ext-iso} and $(**)$.
It follows from $X\in \add(\bigtriangleup_A^c)$ that $F(X)\simeq F'(X)$ and therefore $F(X)\oplus F'(X')\simeq W$ as $B$-modules.
Hence $(\diamondsuit)$ can be written as
$$ 0\lra F(X)\oplus F'(X')\lra F(Y\oplus \bar{Y} \oplus I)\oplus Q\oplus Q'\lraf{g'}F(Z)\lra 0.$$

To complete the proof, we have to show that $F(\bar{Y}\oplus I)\oplus Q \simeq F'(\bar{Y}\oplus P)\oplus J$ for some $J\in \add(\mathscr{I}(B)_{\mathscr{P}})$. Indeed, it follows from $Y\in \add(\bigtriangleup_A^c)$ that $F(Y)\simeq F'(Y)$ as $B$-modules and that both $F(Y)$ and $F'(Y)$ lie in $B\modcat_{\mathscr{I}}$. Since $\bar{Y}$ belongs to $\add(\mathfrak{n}_{F}(A))$ and $I$ belongs to $\add(\mathscr{I}(A)_{\mathscr{P}})$, it follows from Lemma \ref{one-to-one} that $F(\bar{Y} \oplus I)$ lies in $\add(\bigtriangledown_B)$.
Thus $F(\bar{Y} \oplus I)\simeq F(\bar{Y})\oplus F(I) = V\oplus J$ for some $V\in \add(\mathfrak{n}_{G}(B))$ and $J\in \mathscr{I}(B)_{\mathscr{P}}$. Therefore
we have the isomorphisms in $(\stmodcat{B})\modcat$:
$$\begin{array}{rl} \Ext_B^1(-,F(Y\oplus \bar{Y} \oplus I)\oplus Q) & = \; \Ext_B^1(-,F(Y)\oplus F(\bar{Y} \oplus I)\oplus Q)\\ & \simeq \; \Ext_B^1(-,F(Y)\oplus V\oplus Q) \\ & \simeq \; \Ext_B^1(-,F'(Y \oplus\bar{Y} \oplus P)) \quad (\mbox{ by } (\ddag) ).\end{array}$$
As $\bar{Y}\in \add(\mathfrak{n}_{F}(A))$ and $P\in\add(\mathscr{P}(A)_{\mathscr{I}})$, it follows from Lemma \ref{one-to-one} that $F'(\bar{Y} \oplus P)$ is in $\add(\bigtriangleup_B)$ and $F'(\bar{Y} \oplus P)$ is in $B\modcat_{\mathscr{I}}$. Now, Lemma \ref{ext-iso} shows that $F'(Y \oplus\bar{Y}\oplus P)\simeq F(Y)\oplus V\oplus Q$ and $F'(\bar{Y})\oplus F'(P)\simeq V\oplus Q$ as $B$-modules. Thus $F(\bar{Y}) \oplus F(I)\oplus Q \simeq V\oplus J\oplus Q\simeq F'(\bar{Y})\oplus F'(P)\oplus J$ as $B$-modules.
$\square$

\medskip
The following special case of Lemma \ref{exact-seq} is often used in our proofs.
\begin{Koro}\label{exact}
Let $0\ra X\oplus X'\ra  P'\stackrel{g}{\ra}Z\ra 0$
be a minimal exact sequence in $A\modcat$ such that $X\in \add(\bigtriangleup_A^c)$, $X'\in \add(\bigtriangleup_A)$, $P'\in \PI{A}$ and $Z\in A\modcat_{\mathscr{P}}$.
Then there exists a minimal exact sequence
$$0\lra F(X)\oplus F'(X')\lra Q'\lraf{g'}F(Z)\lra 0$$
of $B$-modules with $Q'\in \PI{B}$.
\end{Koro}

\subsection{Stable equivalences induced by almost $\nu$-stable derived equivalences }

In this subsection, we discuss special stable equivalences, called stable equivalences of Morita type.

\begin{Def}\label{def-st-m}
Let $A$ and $B$ be arbitrary finite-dimensional algebras over a field $k$.

$(1)$ $A$ and $B$ are \emph{stably equivalent of Morita type} (see {\rm \cite{broue}}) if there exist bimodules  $_AM_B$ and $_BN_A$ such that

{\rm(i)} $M$ and $N$ are projective as one-sided modules,

{\rm(ii)} $M\otimes_B N\simeq A\oplus P$ as $A$-$A$-bimodules for some projective $A$-$A$-bimodule $P$, and $N\otimes_A M\simeq B\oplus Q$
as $B$-$B$-bimodules for some projective $B$-$B$-bimodule $Q$.

$(2)$ $A$ and $B$ are \emph{stably equivalent of adjoint type} (see {\rm \cite{Xi2008}}) if the bimodules $M$ and $N$ in $(1)$ provide additionally two adjoint pairs $(M\otimes_B-, N\otimes_A-)$ and $(N\otimes_A-, M\otimes_B-)$ of functors on module categories.
\end{Def}

By definition, the module $N$ induces a stable equivalence $N\otimes_A-: \stmodcat{A} \ra \stmodcat{B}$.
Typical examples of stable equivalences of Morita type arise from derived equivalences between self-injective algebras \cite{rickard89}. In this case, derived equivalences induce stable equivalences of Morita type. A generalization of this fact is the class of almost $\nu$-stable derived equivalences introduced in \cite{hx10}. For general Morita theory of derived equivalences and the related notion of tilting complexes, we refer to \cite{Rickard1}.

\begin{Def}{\rm \cite{hx10}} A derived equivalence $F$ of bounded derived module categories between arbitrary Artin algebras $A$ and $B$ with a quasi-inverse $G$ is said to be \emph{almost $\nu$-stable} if the associated radical tilting complexes $\cpx{T}$ over $A$ to $F$ and $\cpx{\bar{T}}$ over $B$ to $G$ are of the form
$$ \cpx{T}: 0\lra T^{-n}\lra \cdots T^{-1}\lra T^0\lra 0 \; \mbox{ and } \,
\cpx{\bar{T}}: 0\lra \bar{T}^{0}\lra \bar{T}^1\lra \cdots \lra \bar{T}^{n}\lra 0,$$
respectively, such that $\add(\bigoplus_{i=1}^n T^{-i})=\add(\bigoplus_{i=1}^n \nu_A( T^{-i}))$ and $\add(\bigoplus_{i=1}^n \bar{T}^{i}) = \add(\bigoplus_{i=1}^n \nu_B(\bar{T}^{i})).$ \end{Def}

It was shown that each almost $\nu$-stable derived equivalence between finite-dimensional algebras over fields induces a stable equivalence of Morita type (see \cite[Theorem 1.1(2)]{hx10}). To get almost $\nu$-stable derived equivalences, we mention the following.

\begin{Lem} {\rm \cite[Corollary 3.14]{hx13}} \label{self-inj}
If $A$ is a finite-dimensional, self-injective algebra over a field and $M\in A\modcat$, then $\End_A(A\oplus M)$ and $\End_A(A\oplus \Omega_A(M))$ are almost $\nu$-stable derived equivalent.
\end{Lem}

Now, we point out that stable equivalences of adjoint type have some nice properties.
\begin{Lem}\label{lemast} Let $A$ and $B$ be arbitrary finite-dimensional algebras over a field.
Suppose $A$ and $B$ are stably equivalent of adjoint type induced by $_AM_B$ and $_BN_A$. Write $_AM\otimes_BN_A\simeq A\oplus P$ and $_BN\otimes_AM_B\simeq B\oplus Q$ as bimodules. Then the following hold:

$(1)$ $\add(\nu_AP)=\add(_AP)$ and $\add(\nu_BQ)=\add(_BQ)$, where $\nu_A$ is the Nakayama functor of $A$.

$(2)$ If $S$ is a simple $A$-module with $\Hom_{A}(P,S)=0$, then $N\otimes_AS$ is a simple $B$-module with $\Hom_B(Q, N\otimes_AS)=0$.

$(3)$ For an $A$-module $X$ and a $B$-module $Y$, we have $\Omega_B^i(N\otimes_AX)\simeq N\otimes_A\Omega_A^i(X)$ in $\stmodcat{B}$ and $\Omega_A^i(M\otimes_BY)\simeq M\otimes_B\Omega_B^i(Y)$ in $\stmodcat{A}\,$ for $i\ge 0$.

$(4)$ For an $A$-module $X$ and a $B$-module $Y$, there are equalities $$\del(_AM\otimes_BN\otimes_AX)=\del({_B}N\otimes_AX)=\del(_AX)\mbox{ and } \del(_BN\otimes_AM\otimes_BY)=\del({_A}M\otimes_BY)=\del(_BY).$$
\end{Lem}
{\it Proof.} (1) and (2) follow from \cite[Lemma 3.1]{hx18}, while (3) and (4) can be deduced easily. $\square$

\begin{Lem}\label{seat}
Let $A$ and $B$ be arbitrary Artin algebras with no separable direct summands. If $A$ and $B$ are stably equivalent of Morita type, then they are even stably equivalent of adjoint type.
\end{Lem}
{\it Proof.}
Since $A$ and $B$ are stably equivalent of Morita type and have no separable direct summands, it follows from \cite[Proposition 2.1 and Theorem 2.2]{liu08} which are valid also for Artin algebras by checking the argument there, that $A$ and $B$ have the same number of indecomposable direct summands (as two-sided ideals) and that we may write $A=A_1\times A_2\times\cdots \times A_s$ and $B = B_1\times B_2 \times \cdots \times B_s$ as products of indecomposable algebras, such that the blocks $A_i$ and $B_i$
are stably equivalent of Morita type for $1\le i\le s$. Suppose that $M^{(i)}$ and $N^{(i)}$ define a stable equivalence of Morita type between $A_i$ and $B_i$. Observe that the two results \cite[Lemma 2.1 and Corollary 3.1]{D-MV2007} hold true for indecomposable, non-separable Artin algebras. Thus, by \cite[Lemma 2.1]{D-MV2007}, we may assume that $_{A_i}{M^{(i)}}_{B_i}$ and $_{B_i}{N^{(i)}}_{A_i}$ are indecomposable, non-projective bimodules. Since the algebra $A_i$ is indecomposable and non-separable by assumption, it follows from \cite[Corollary 3.1]{D-MV2007} that $(M^{(i)}\otimes_{B_i}-, N^{(i)}\otimes_{A_i}-)$ and $(N^{(i)}\otimes_{A_i}-, M^{(i)}\otimes_{B_i}-)$ are adjoint pairs between $A_i\modcat$ and $B_i\modcat$. Let $M:=\bigoplus_{1\le j\le s}M^{(j)}$ and $N:=\bigoplus_{1\le j\le s}N^{(j)}$.
Then $_AM_B$ and $_BN_A$ define a stable equivalence of adjoint type between $A$ and $B$.
$\square$

\subsection{Elimination of nodes in algebras}
In this section, we recall a general procedure of eliminating nodes of algebras introduced by Martinez-Villa \cite{MV1980}, and prove some new properties of this elimination of nodes.

Recall that a non-projective, non-injective simple $A$-module $M$ is called a \emph{node} of $A$ if the middle term $P$ of an almost split sequence starting at $M$, $0\ra M\ra P\ra {\rm Tr\,}D(M)\ra 0$, is projective. Note that $A$ has no nodes if and only if $A^{\opp}$ has no nodes \cite[Lemma 1]{MV1980}.

In \cite[Theorem 2.10]{MV1980}, Mart\'inez-Villa showed that any Artin algebra with nodes is stably equivalent to an Artin algebra without nodes. The process of removing nodes runs precisely as follows. Suppose that $A$ is an Artin algebra with nodes. Let $\{S(1), S(2), \cdots, S(n)\}$ be a complete set of non-isomorphic simple $A$-modules. Suppose that $P(i)=Ae_i$ has the top $S(i)$ with $e_i^2=e_i\in A$ for $1\le i\le n$. We may assume that $\{S(1),\cdots,S(m)\}$ is a complete set of nodes of $A$ with $m\le n$. Set $S:=\bigoplus_{i=1}^mS(i)$. Let $I$ be the trace of $S$ in $A$. Then $I\in \add(_AS)$. By the definition of nodes, $S\in \add(\soc(_AA))$ and $S\in \add(_AI)$. Thus $\add(_AI)=\add(_AS)$. Clearly, $I^2=0$ and $\rad(A)I=0$. Let $J:=$ann$_l(I)$ be the left annihilator of $I$. Then $\rad(A)\subseteq J$ and $A/J$ is semisimple. Since $_AI$ has only composition factors $S(i)$ for $1\le i\le m$, we have $e_iI\neq 0$ for $1\le i\le m$ and $e_jI=0$ for $m+1\le j\le n$. This yields $\add(_AA/J)=\add(_AS)$. Note that $I$ is a two-sided ideal of $A$ and $JI=0$. Thus $I$ is an $(A/J)$-$(A/I)$-bimodule. We may form the triangular matrix algebra
$$A':=\begin{pmatrix}
A/I&0\\
I&A/J
\end{pmatrix}.$$
The algebra $A'$ is called the \emph{model algebra} of $A$. The following lemma describes some common properties of $A$ and $A'$.

\begin{Lem}\label{node}
$(1)$ The triangular matrix Artin algebra $A'$ has no nodes.

$(2)$ $A$ and $A'$ are stably equivalent.

$(3)$ $A$ and $A'$ have the same number of non-isomorphic, non-projective simples.

$(4)$ If $A$ is Frobenius-finite, then so is $A'$.
\end{Lem}

{\it Proof.} The first two statements are taken from \cite[Theorem 2.10]{MV1980}, while we prove (3) and (4).

(3) As is known, $A'$-modules can be identified with triples $(X,Y,f)$, where $X$ is an $A/I$-module, $Y$ is an $A/J$-module and $f:I\otimes_{A/I}X\ra Y$ is a homomorphism of $A/J$-modules. It follows from $I^2=0$ that $I\subseteq \rad(A)$. Thus simple $A$-modules coincide with simple $A/I$-modules, and therefore $A$ and $A/I$ have the same number of non-isomorphic simple modules. Note that the projective cover of an $A/I$-module $X$ is of the form $P/IP$ with $P$ being a projective cover of the $A$-module $_AX$. Obviously, the simple $A'$-modules are either of the form $(T,0,0)$, where $T$ is a simple $A$-module, or of the form $(0,T',0)$, where $T'$ is a simple $A/J$-module. The indecomposable projective $A'$-modules are either of the form $\tilde{P}:=(P/IP,I\otimes_{A/I}P/IP, {\rm id})$ with $P$ an indecomposable projective $A$-module or of the form $(0,T',0)$ with $T'$ an indecomposable projective $A/J$-module. Thus $(0,T',0)$ is a projective simple $A'$-module, and so the indecomposable non-projective simple $A'$-modules are of the form $(T,0,0)$, where $T$ is a simple $A$-module.

We prove that $A$ and $A'$ have the same number of non-isomorphic, non-projective simples. Indeed, take a simple $A$-module $T$, then $IT=0$ and $(T,0,0)$ is a simple $A'$-module. If $T$ is a projective $A$-module, then $T$ is also a projective $A/I$-module, and therefore $(T,0,0)$ is a projective $A'$-module. Thus $(T,0,0)$ is a projective simple $A'$-module. Suppose that $T$ is not a projective $A$-module. Let $P(T)$ be a projective cover of $_AT$. Then $P(T)\not\simeq T$. If $IP(T)=0$, then $P(T)$ is a projective cover of $_{A/I}T$. Thus $(P(T),0,0)$ is a projective cover of $(T,0,0)$ and $(T,0,0)$ is not a projective $A'$-module. If $IP(T)\neq 0$, then it follows
that $I\otimes_{A/I}\big(P(T)/IP(T)\big)\simeq I\otimes_A\big(P(T)/IP(T)\big)\simeq I\otimes_A(A/I)\otimes_AP(T)\simeq I\otimes_AP(T)\simeq IP(T)\neq 0$. Thus the $A'$-module $\big(P(T)/IP(T),I\otimes_{A/I}P(T)/IP(T),{\rm id}\big)$ is a projective cover of $(T,0,0)$. This implies that $(T,0,0)$ is not projective. Thus $A$ and $A'$ have the same number of non-isomorphic, non-projective simples.

(4) Suppose that $(_{A/I}X, _{A/J}Y,f)$ is an indecomposable $A'$-module in $A'\stp$. We show that $_{A/I}X\in A/I\stp$ and $_{A/J}Y\in A/J\stp$. Indeed, $(_{A/I}X, _{A/J}Y,f)$ is projective-injective with $\nu_{A'}(_{A/I}X, _{A/J}Y,f)\in A'\stp$. It follows from \cite[Proposition 2.5, p.76]{Aus1997} that there are two possibilities:

(a) $_{A/J}Y=0$ and $_{A/I}X$ is an indecomposable projective-injective $A/I$-module with $I\otimes_{A/I}X=0$;

(b) $_{A/I}X=0$ and $_{A/J}Y$ is an indecomposable projective-injective $A/J$-module with $\Hom_{A/J}(I,Y)=0$.

Suppose (a) holds. Let $T_0:=\top(_{A/I}X)$. Then $\nu_{A/I}(X)$ is an injective envelope of $_{A/I}T_0$. By \cite[Proposition 2.5, p.76]{Aus1997}, we see that $_{A'}(T_0,0,0)$ is a simple $A'$-module, that $_{A'}(X,0,0)$ is a projective cover of $_{A'}(T_0,0,0)$, and that $_{A'}(\nu_{A/I}(X),0,0)$ is an injective envelope of $_{A'}(T_0,0,0)$. Thus $\nu_{A'}(_{A/I}X, 0,0)\simeq (\nu_{A/I}(X),0,0)\in A'\stp$. This implies that $\nu^i_{A/I}(X)$ is projective-injective
for all $i\ge 0$, and $X\in A/I\stp$. Similarly, if (b) holds, then $Y\in A/J\stp$.

Let $S$ be a direct sum of representatives of the isomorphism classes of nodes of $A$. Since $I$ is the trace of $S$ in $A$ and $J$ is the left annihilator of $I$, we have $\add(_AI)=\add(_AS)=\add(_AA/J)$. Thus $\Hom_A(I,Z)=\Hom_{A/J}(I,Z)\neq 0$ for any $A/J$-module $Z\ne 0$. Hence we can assume that $\{(X_1,0,0),\cdots (X_r,0,0)\}$ is a complete set of all non-isomorphic indecomposable modules in $A'\stp$ for some natural number $r$ and $X_i\in A/I\stp$ with $I\otimes_{A/I}X_i=0$ for $1\le i\le r$. Since each indecomposable projective $A/I$-module is of the form $P/IP$ for some indecomposable projective $A$-module $P$, we assume $X_i\simeq P_i/IP_i$ for some indecomposable projective $A$-module $P_i$. Note that $I\otimes_{A/I}\big(Q/IQ\big)\simeq I\otimes_A\big(Q/IQ\big)\simeq I\otimes_A(A/I)\otimes_AQ\simeq I\otimes_AQ\simeq IQ$ for each projective $A$-module $Q$. It follows from $I\otimes_{A/I}X_i=0$ that $IP_i=0$, and therefore $X_i\simeq P_i$ as $A/I$-modules and $\soc(_AP_i)$ is a simple $A$-module. Since $P_i$ is a projective $A$-module, the trace of $S$ in $P_i$ is equal to $IP_i$. It then follows from $IP_i=0$ that $\soc(_AP_i)$ has no nodes as its direct summands for $1\le i\le r$. Set $U:=\bigoplus_{i=1}^rP_i$. Then $\soc(_AU)$ has no nodes as its direct summands. Since $\nu_{A'}(U,0,0)$ is $\nu$-stably projective and $\nu_{A'}(U,0,0)\simeq (\nu_{A/I}(U),0,0)$,  there are isomorphisms: $U\simeq \nu_{A/I}(U)$ and $\top(_{A/I}U)\simeq\soc(_{A/I}U)$ as $A/I$-modules. Thus $\top(_{A}U)$ is isomorphic to $\soc(_{A}U)$ and has no nodes as its direct summands. In particular, $\Hom_A(U,I)=0$.
Applying $\Hom_A(U,-)$ to the exact sequence
$$0\lra I\lra A\lra A/I\lra 0$$
of $A$-$A$-bimodules, we get the  exact sequence of $A^{\opp}$-modules
$$0\lra \Hom_A(U,I)\lra \Hom_A(U,A)\lra \Hom_A(U,A/I)\lra 0.$$
It follows from $\Hom_A(U,I)=0$ that $\Hom_A(U,A)\simeq \Hom_A(U,A/I)$. Clearly, $\Hom_A(U,A/I)= \Hom_{A/I}(U,A/I)$ as $A^{\opp}$-modules. Thus $D\Hom_A(U,A)\simeq D\Hom_{A/I}(U,A/I)$ as $A$-modules. As $\nu_{A/I}(U)=D\Hom_{A/I}(U,A/I)\simeq U$ as $A$-modules and $D\Hom_A(U,A)\simeq U$ as $A$-modules, we get $U\in A\stp$. Let $\Lambda$ be the Frobenius part of $A$, and let $\Lambda'$ be the Frobenius part of $A'$. Then $\End_A(U)$ is of the form $f\Lambda f$ for an idempotent $f\in \Lambda$, and $$\Lambda':=\End_{A'}((U,0,0))\simeq \End_{A/I}(U)\simeq \End_{A}(U).$$
If $A$ is Frobenius-finite, then $\End_A(U)$ is representation-finite, and therefore $A'$ is Frobenius-finite.
$\square$

\section{New invariants of stable equivalences}\label{proof-main}
This section is devoted to the proofs of Theorems \ref{arc+del}-\ref{main-thm}.
We keep the notation introduced in the previous sections.

Let $A$ be an Artin algebra over a commutative Artin ring $k$. Recall that a projective $A$-module $P$ is said to be \emph{$\nu$-stably projective} if $\nu_A^i(P)$ is projective for all $i> 0$. Here $\nu_A$ is the Nakayama functor $D\Hom_A(-,A)$ of $A$. Let $U$ be a direct sum of representatives of the isomorphism classes of indecomposable $\nu$-stably projective $A$-modules. Since $\nu_A(U)$ is $\nu$-stably projective, we have $U\simeq \nu_A(U)$ and $\top(U)\simeq\soc(U)$, where $\soc(U)$ is the socle of the $A$-module $U$. Clearly,  $\soc(U)\simeq \Omega_A(U/\soc(U))\oplus Q$ for some projective $A$-module $Q$. Thus $\del(\top(U))=\del(\soc(U))=0$ by definition. Let $V$ be a direct sum of representatives of the isomorphism classes of
indecomposable projective $A$-modules that are neither simple nor $\nu$-stably projective. Then $\del(A)=\del(\top(U\oplus V))
=\max\{\del(\top(U)), \del(\top(V))\}
=\del(\top(V))$.

\subsection{Homological data under stable equivalences \label{sect4.1}}

Let $A$ and $B$ be Artin $k$-algebras that have neither nodes nor semisimple direct summands. Assume that $F: \stmodcat{A}\ra \stmodcat{B}$ is an equivalence of $k$-categories.

Under these assumptions, $\mathfrak{n}_{F}(A)=\varnothing$, $\mathfrak{n}_{F^{-1}}(B)=\varnothing$, where $\mathfrak{n}_{F}(A)$ denotes the set of isomorphism classes of $F$-exceptional nodes of $A$ (see Subsection \ref{sect3.1}), and there is a bijection $F': \mathscr{P}(A)_{\mathscr{I}}\ra \mathscr{P}(B)_{\mathscr{I}}$ (see  Lemma \ref{one-to-one}). Further, Lemma \ref{exact-seq} can be specified as follows.
\begin{Lem}\label{seq}{\rm\cite[Theorem 1.7]{MV1990}} Let $0\ra X\oplus P_1\lraf{f} Y\oplus P\oplus P'\lraf{g} Z\ra 0$ be a minimal exact sequence of $A$-modules, where $X,Y,Z\in A\modcat_{\mathscr{P}}$, $P_1,P\in \mathscr{P}(A)_{\mathscr{I}}$  and $P'$ is a projective-injective $A$-module. Then
there is a minimal exact sequence
$$0\lra F(X)\oplus F'(P_1)\lraf{f'} F(Y)\oplus F'(P)\oplus Q\lraf{g'}F(Z)\lra 0$$
in $B\modcat$ with $Q$ a projective-injective $B$-module and $\underline{g'} =F(\underline{g})$. In particular, $\Omega_{B}F(Z)\simeq F\Omega_A(Z)$ in $\stmodcat{B}$ for $Z\in A\modcat_{\mathscr{P}}$.
\end{Lem}

\begin{Lem}\label{delalp}
For $X\in A\modcat_{\mathscr{P}}$, we have $\del(X)=\del(F(X))$.
\end{Lem}
{\it Proof.} We show $\del(F(X))\le \del(X)$. We may assume $d:=\del(X)<\infty$. Then, by the definition of delooping levels, there exists $M\in A\modcat_{\mathscr{P}}$ such that $\Omega_A^d(X)\in \add(_AA\oplus \Omega_A^{d+1}(M))$. Thus $F(\Omega_A^d(X))\in \add({_{B}}B\oplus F\Omega_A^{d+1}(M))$ by the additivity of the functor $F$. On the other hand, it follows from Lemma \ref{seq} that $\Omega_{B}F(X) \simeq F\Omega_A(X)$ and $\Omega_{B}F(M)$ $\simeq F\Omega_A(M)$ in $\stmodcat{B}$. Then $\Omega_{B}^i F(X)\simeq F\Omega_A^i(X)$ and $\Omega_{B}^i F(M)\simeq F\Omega_A^i(M)$ in $\stmodcat{B}$ for $i\ge 0$. Thus  $\Omega_{B}^d F(X)\in \add({_{B}}B\oplus\Omega_B^{d+1}F(M))$, and therefore $\del(F(X))\le d=\del(X)<\infty$. Similarly, we show $\del(X)\le\del(F(X)).$ Thus $\del(X)=\del(F(X))$. $\square$

\medskip
{\bf Proof of Theorem \ref{arc+del}}. Suppose that $A$ and $B$ are stably equivalent Artin algebras without nodes.

(i) $\del(A)=\del(B).$ In fact, the delooping levels of algebras involve only simple modules. If $A$ or $B$ has a semisimple direct summand, then the simple modules belonging to the semisimple direct summand have delooping levels $0$, and therefore do not contribute to the delooping levels of the algebra considered. So we may remove all semisimple direct summands from both algebras $A$ and $B$. Of course, the resulting algebras are still stably equivalent. Thus we assume that both $A$ and $B$ do not have any semisimple direct summands. Let $V$ be a direct sum of representatives of the isomorphism classes of indecomposable projective $A$-modules that are neither simple nor $\nu$-stably projective, and let $V'$ be a direct sum of representatives of the isomorphism classes of indecomposable projective $B$-modules that are neither simple nor $\nu$-stably projective. It follows from \cite[Lemma 2.5]{MV1990}, which is true also for Artin algebras, that $F(\top(V))\simeq\top(V')$ as $B$-modules. Note that $\top(V)$ does not have any nonzero projective direct summands. By Lemma \ref{delalp},  we have $\del(\top(V))=\del(\top(V'))$. Thus $\del(A)=\del(\top(V))
=\del(\top(V'))=\del(B)$.

(ii) $\phi\dim(A)=\phi\dim(B)$ and $\psi\dim(A)=\psi\dim(B)$. For Artin algebras $A_1$ and $A_2$, there are equalities $\phi\dim(A_1\times A_2)=\max\{\phi\dim(A_1), \phi\dim(A_2)\}$ and $\psi\dim(A_1\times A_2) = \max\{\psi\dim(A_1), \psi\dim(A_2)\}$. Since $\phi$- and $\psi$-dimensions of semisimple algebras are $0$, we may remove semisimple direct summands from $A$ and $B$ if they have any. Then the resulting algebras are still stably equivalent. So we may assume that both algebras $A$ and $B$ do not have any semisimple direct summands. Let $F: \stmodcat{A}\ra \stmodcat{B}$ be a stable equivalence between $A$ and $B$. We also denote by $F$ the correspondence from $A$-modules to $B$-modules, which takes projective $A$-modules to $0$. As a functor of $k$-categories, $F$ is additive and commutes with finite direct sums in $A\stmodcat$. Thus the map $\tilde{F}: K(A)\ra K(B)$ given by $\tilde{F}([X]):=[F(X)]$, is a well-defined homomorphism of Grothendieck groups. It is actually an isomorphism of abelian groups. By Lemma \ref{seq} we have $\Omega_{B}(F(X))\simeq F(\Omega_A(X))$ in $B\stmodcat$ for $X\in A\modcat$. Let $\langle X\rangle$ be the $\mathbb{Z}$-submodule of $K(A)$ generated by the isomorphism classes of indecomposable, non-projective direct summands of $X$. For $n\ge 0$, the following diagrams are commutative
$$\xymatrix{
K(A)\ar[r]^-{\tilde{F}}\ar[d]^-{\Omega_A}
& K(B)\ar[d]^-{\Omega_{B}}&
\Omega_A^n\langle X\rangle\ar[r]^-{\tilde{F}_{res}}\ar[d]^-{\Omega_A}
&\Omega_{B}^n\langle F(X)\rangle\ar[d]^-{\Omega_{B}}\\
K(A)\ar[r]^-{\tilde{F}}&K(B),&
\Omega_A^{n+1}\langle X\rangle\ar[r]^-{\tilde{F}_{res}}&\Omega_{B}^{n+1}\langle F(X)\rangle
}$$
where $\tilde{F}_{res}$ is the restriction of $\tilde{F}$. Since $\tilde{F}:K(A)\ra K(B)$ is an isomorphism of abelian groups, the $\mathbb{Z}$-module homomorphism $\Omega_A:\Omega_A^n\langle X\rangle\ra \Omega_A^{n+1}\langle X\rangle$ is an isomorphism for $n\ge 0$ if and only if so is the $\mathbb{Z}$-homomorphism $\Omega_{B}:\Omega_{B}^n\langle F(X)\rangle\ra \Omega_{B}^{n+1}\langle F(X)\rangle$ for $n\ge 0$. By the definition of $\phi$-dimensions, $\phi(X)=\phi(F(X))$ and $\phi\dim(A)\le \phi\dim(B)$. Similarly, $\phi\dim(B)\le \phi\dim(A)$. Thus $\phi\dim(A)=\phi\dim(B)$. For $Y\in A\modcat$, since $\Omega_{B}(F(Y))\simeq F(\Omega_A(Y))$ in $B\stmodcat$ and $F$ is an equivalence, we get $\pd(_{B}F(Y))=\pd(_AY)$. Then $\psi(X)=\psi(F(X))$ and $\psi\dim(A)=\psi\dim(B)$. $\square$

\medskip
Theorem \ref{arc+del} may fail if Artin algebras have nodes. This can be seen by the following examples.

\begin{Bsp}{\rm\label{exnode}
(1) Let $A_1$ be the algebra over a field $k$, given by the quiver with a relation:
$$\xymatrix{
1\,\bullet\ar@(ru,rd)^{\alpha \; ,}&\quad \alpha^2=0.
}$$
Clearly, $A_1$ has a node and is stably equivalent to the path algebra $A'_1$ given by the quiver $1\bullet\leftarrow \bullet 2$.
Note that $A'_1$ has no nodes and its Frobenius part is $0$.  In this case, both $A_1$ and $A'_1$ have only $1$ non-projective simple module. Clearly, $\del(A_1)=\phi\dim(A_1)=\psi\dim(A_1)=0 < 1=\del(A'_1)=\phi\dim(A'_1)=\psi\dim(A'_1).$ Remark that $A_1$ and $A'_1$ are never stably equivalent of Morita type by Lemma \ref{thmmt} below.

(2) Let $A_2$ be the algebra over a field $k$, given by the quiver with a relation:
$$\xymatrix{
A_2: \quad 1\,\bullet \ar@<0.5ex>[r]^{\alpha} & \bullet 2\;,\ar@<0.5ex>[l]^{\beta}&\alpha\beta=0
}$$
(see \cite[Example 4.9]{Xi2021} for general situations).We consider the  $2$ almost split sequences in $A_2\modcat$
$$0\lra S(1)\lra P(2)\lra S(2)\lra 0 \; \mbox{ and  } \; 0\lra S(2)\lra I(2)\lra S(1)\lra 0,$$ where $P(i)$, $I(i)$ and $S(i)$ are the projective, injective and simple modules corresponding to the vertex $i$, respectively. Clearly, $A_2$ has the Frobenius part isomorphic to $A_1$, and a unique node $S(1)$. Let $I$ be the trace of $S(1)$ in $A_2$ and $J$ be the left annihilator of $I$ in $A_2$. Then $I=\{r_1\beta\alpha+r_2\beta\mid r_1,r_2\in k\}$ and $J=\{r_1e_2+r_2\alpha+r_3\beta\mid r_i\in k, 1\le i\le 3\}$. Define $A'_2$ to be the triangular matrix algebra
$$A'_2=\begin{pmatrix}
A_2/I&0\\
I&A_2/J
\end{pmatrix}\simeq
\begin{pmatrix}
k&0&0\\
k&k&0\\
k&k&k
\end{pmatrix}.$$
Then $A'_2$ has no nodes and its Frobenius part is $0$. By Lemma \ref{node}(2), $A_2$ and $A'_2$ are stably equivalent. Particularly, they have $2$ non-isomorphic, non-projective simple modules, and $\del(A_2) = \phi\dim(A_2)=\psi\dim(A_2)=2>1$ = $\del(A'_2)= \phi\dim(A'_2) = \psi\dim(A'_2).$
}\end{Bsp}

\medskip
For stable equivalences of Morita type, the requirement that algebras considered have no nodes can be eliminated from Theorem \ref{arc+del}.

\begin{Lem}\label{thmmt} Let $A$ and $B$ be arbitrary finite-dimensional algebras over a field $k$. If $A$ and $B$ are stably equivalent of Morita type, then

$(1)$ $\del(A)=\del(B)$.

$(2)$ $\phi\dim(A)=\phi\dim(B)$ and $\psi\dim(A)=\psi\dim(B).$
\end{Lem}

{\it Proof.}
(1) Let $A=\Lambda_1\times \Lambda_2$ and $B=\Gamma_1\oplus\Gamma_2$, where $\Lambda_1$ and $\Gamma_1$ are separable algebras, and where $\Lambda_2$ and $\Gamma_2$ are algebras without separable direct summands. Since $A$ and $B$ are stably equivalent of Morita type, it follows from the proof of \cite[Theorem 4.7]{lx05} that $\Lambda_2$ and $\Gamma_2$ are stably equivalent of Morita type. By \cite[Proposition 2.1]{liu08}, we can suppose $\Lambda_2=A_1\times \cdots \times A_s$ and $\Gamma_2=B_1\times \cdots \times B_s$, where all $A_i$ and $B_i$ are indecomposable algebras. By \cite[Theorem 2.2]{liu08},  we may assume that $A_i$ and $B_i$ are stably equivalent of Morita type for all $i$ (up to re-ordering the summands). Since $A_i$ and $B_i$ are stably equivalent, $A_i$ is a semisimple algebra if and only if so is $B_i$ for $1\le i\le s$. Thus we may suppose that $A_1, \cdots, A_t, B_1,\cdots, B_t$ are non-semisimple algebras and that $A_{t+1}, \cdots, A_s, B_{t+1}, \cdots, B_s$ are semisimple algebras. Let $A_0=A_1\times \cdots\times A_t$ and $B_0=B_1\times \cdots \times B_t$. Then $A_0$ and $B_0$  do not have semisimple direct summands. Since $A_i$ and $B_i$ are stably equivalent of Morita type for $1\le i\le t$, we see that $A_0$ and $B_0$ are stably equivalent of Morita type. Suppose that two bimodules $_{A_0}M_{B_0}$ and $_{B_0}N_{A_0}$ define a stable equivalence of Morita type between $A_0$ and $B_0$. Then ${_{A_0}}M\otimes_{B_0}N_{A_0}\simeq A_0\oplus P$ and $_{B_0}N\otimes_{A_0}M_{B_0}\simeq B_0\oplus Q$, where $_{A_0}P_{B_0}$ and $_{B_0}Q_{B_0}$ are projective bimodules. We may assume that $M$ and $N$ have no nonzero projective summands as bimodules by \cite[Lemma 4.8(1)]{lx05}. Thus $M$ and $N$ define a stable equivalence of adjoint type between $A_0$ and $B_0$ by \cite[Lemma 4.1(2)]{cpx1}. It then follows from \cite[Lemma 3.1(2)]{hx18} that $\add(\nu_{A_0}P)=\add(_{A_0}P)$. Thus $\del(\top(P))=\del(\soc(P))=0$ and $$\del(A_0)=\max\{\del(_{A_0}S)\mid {_{A_0}}S \mbox{ is simple with } \Hom_{A_0}(P,S)=0\}.$$
Let $S$ be a simple $A_0$-module with $\Hom_{A_0}(P,S)=0$. By Lemma \cite[Lemma 3.1(5)]{hx18}, $_{B_0}N\otimes_{A_0}S$ is a simple $B_0$-module. It follows from Lemma \ref{lemast}(4) that $\del(_{A_0}S)=\del(_{B_0}N\otimes_{A_0}S)\le \del(B_0)$, and therefore $\del(A_0)\le \del(B_0)$. Similarly, we can show $\del(B_0)\le \del(A_0)$. Thus $\del(A_0)=\del(B_0)$. Since the delooping levels of semisimple blocks are $0$, we have $\del(A)=\del(A_0)=\del(B_0)=\del(B)$.

(2) Suppose that $A$ and $B$ are stably equivalent of Morita type defined by $_AM_B$ and $_BN_A$. Since $_BN$ is projective, the functor $F:=N\otimes_A-:A\modcat\ra B\modcat$ takes projective $A$-modules to projective $B$-modules, and commutes with finite direct sums. Thus $F$ induces an equivalence: $\stmodcat{A}  \ra B\stmodcat$. As in the proof of Theorem \ref{arc+del}(ii), we obtain $\phi\dim(A)=\phi\dim(B)$ and $\psi\dim(A)=\psi\dim(B).$
$\square$

\medskip
\begin{Rem}\label{rmk-tilt} {\rm In contrast to stable equivalences of Morita type, tilting procedure preserves neither delooping levels nor $\phi$-dimensions.
For example, the path algebra $A$ (over a field) of the quiver $\bullet\lraf{\alpha} \bullet\lraf{\beta} \bullet$ can be tilted to the quotient algebra $B:=A/(\alpha\beta)$. In this case, $A$ has no nodes, but $B$ has a node, while we have $\del(A)=\phi\dim(A)=\psi\dim(A)=1<2=\del(B)=\phi\dim(B)=\psi\dim(B)$. This shows that in general derived equivalences do not have to preserve the delooping levels and the $\phi$-dimensions of algebras. Nevertheless, we will show that almost $\nu$-stable derived equivalences do preserve delooping levels, $\phi$-dimensions and $\psi$-dimensions in the next section.
}\end{Rem}

\subsection{Homological data under almost $\nu$-stable derived equivalences}
In this subsection, we apply Lemma \ref{thmmt} to show Theorem \ref{a-n-s-d-e}.

\medskip
{\bf Proof of Theorem \ref{a-n-s-d-e}.} $(1)$ Suppose that there is an almost $\nu$-stable derived equivalence between finite-dimensional algebras $\Lambda$ and $\Gamma$ over a field $k$.
It follows from \cite[Theorem 1.1(2)]{hx10} that $\Lambda$ and $\Gamma$ are stably equivalent of Morita type. By Lemma \ref{thmmt}, $\del(\Lambda)=\del(\Gamma)$, $\phi\dim(\Lambda)=\phi\dim(\Gamma)$ and $\psi\dim(\Lambda)=\psi\dim(\Gamma)$.

$(2)$ Let $A$ be a self-injective algebra over a field and $X\in A\modcat$. We show that
$$(a) \; \del\big(\End_A(A\oplus X)\big) = \del\big(\End_A(A\oplus \Omega_A(X))\big)=\del\big(\End_A(A\oplus D{\rm Tr}(X))\big).$$

In fact, by Lemma \ref{self-inj} (see also the remark at the end of Section 3 in \cite{hx10}), we know that $\End_A(A\oplus X)$ and $\End_A(A\oplus \Omega_A(X))$ are almost $\nu$-stable derived equivalent and therefore they are stably equivalent of Morita type. Hence it follows from Lemma \ref{thmmt} that
$$(\sharp)\quad \del\big(\End_A(A\oplus X)\big)=\del\big(\End_A(A\oplus \Omega_A(X))\big).$$As $\nu_A$ is an auto-equivalence of $A\modcat$ and $D{\rm Tr}(Y)\simeq \Omega^2(\nu_A(Y))$ in $A\stmodcat$ for $Y\in A\modcat$, we have
$$\begin{array}{rl} \del\big(\End_A(A\oplus X)\big) = & \del\big(\End_A(\nu_A(A\oplus X))\big) = \del\big(\End_A(A\oplus \nu_A(X))\big)\\ = &\del\big(\End_A(A\oplus \Omega(\nu_A(X))\big)\; \; ( \mbox{ by } (\sharp))\\
 = & \del\big(\End_A(A\oplus \Omega^2(\nu_A(X))\big) \; \; (\mbox{ by } (\sharp))\\
 = &\del\big(\End_A(A\oplus D{\rm Tr}(X))\big),\end{array}$$where the last equality is due to the fact that $\End_A(A\oplus \Omega^2(\nu_A(X))$ and $\End_A(A\oplus D{\rm Tr}(X))$ are Morita equivalent.

Similarly, we can prove the equalities for $\phi$-dimensions and $\psi$-dimensions by Lemma \ref{thmmt}(2):

$(b)\; \phi\dim\big(\End_A(A\oplus X)\big) = \phi\dim\big(\End_A(A\oplus \Omega_A(X))\big)=\phi\dim\big(\End_A(A\oplus D{\rm Tr}(X))\big).$

$(c)\; \psi\dim\big(\End_A(A\oplus X)\big) = \psi\dim\big(\End_A(A\oplus \Omega_A(X))\big)=\psi\dim\big(\End_A(A\oplus D{\rm Tr}(X))\big).$ $\square$

\medskip
Theorem \ref{a-n-s-d-e} can distinguish almost $\nu$-stable derived equivalences among derived equivalences. For instance, let $A$ and $B$ be algebras given by the following quivers $Q_A$ and $Q_B$ with relations, respectively:
$$\begin{array}{ccc}
\xymatrix{
Q_A: \quad \bullet\ar@<2.2pt>[r]^-{\alpha}^(0.27){1}^(1){2}& \bullet\ar@<2.5pt>[l]^-{\delta}\ar@<2.5pt>[r]^-{\beta}^(1){3} &\bullet\ar@<2.5pt>[l]^{\gamma}}
 &\hspace{1cm}& Q_B: \;
\xymatrix{
\bullet\ar[r]^{\alpha'}^(0){1}^(1){2}& \bullet\ar[dl]^{\beta'}\\
\bullet\ar[u]^{\gamma'}_(-0.1){3} }\\
\alpha\delta\alpha=\gamma\delta=\delta\alpha-\beta\gamma=0\; ;
&&
\quad \quad \alpha'\beta'\gamma'\alpha'=\gamma'\alpha'\beta'\gamma'=0. \\ \end{array}$$
It was shown in \cite[Example 4.10]{hx} that $A$ and $B$ are derived equivalent. One can check that both algebras have no nodes and $\del(A)=2\ne 1=\del(B)$. Thus $A$ and $B$ are neither almost $\nu$-stable derived equivalent by Theorem \ref{a-n-s-d-e} nor stably equivalent by Theorem \ref{arc+del}.

\subsection{Stable equivalences of algebras and their Frobenius parts \label{sect3.2}}

Now we turn to the proof of Theorem \ref{main-thm}. We start with the following lemma.

\begin{Lem}\label{bijective} Let $F:\stmodcat{A}\ra \stmodcat{B}$ define a stable equivalence between Artin algebras $A$ and $B$, and let $G$ be a quasi-inverse of $F$.
If $\ndd(A)\ge 1$ and $\ndd(B)\ge 1$, then there exist bijections
$$F:\mathscr{I}(A)_{\mathscr{P}}\lra \mathscr{I}(B)_{\mathscr{P}},\;
F:\mathfrak{n}_{F}(A)\lra \mathfrak{n}_{G}(B),\;
F':\mathscr{P}(A)_{\mathscr{I}}\lra \mathscr{P}(B)_{\mathscr{I}}
\mbox{ and }
F':\mathfrak{n}_{F}(A)\lra \mathfrak{n}_{G}(B).$$
\end{Lem}
{\it Proof.} Suppose $I\in \mathscr{I}(A)_{\mathscr{P}}$, we show $F(I)\in \mathscr{I}(B)_{\mathscr{P}}$. Indeed, let $S$ be the socle of $I$. By Lemma \ref{pi=stp}(2)-(3), $S$ is not injective. Thus $S\not\simeq I$ and the natural projection $\pi:I\ra I/S$ is an irreducible map. 
Since $I$ is not a projective module, we have $I/S\in A\modcat_{\mathscr{P}}$. Thus $0\neq\pi\in \stmodcat{A}$ and $0\neq F(\pi)\in \stmodcat{B}$. By \cite[Chapter X, Proposition 1.3]{Aus1997}, $F(\pi):F(I)\ra F(I/S)$ is irreducible. By Lemma \ref{one-to-one}, we have $F(I)\in\bigtriangledown_A$, namely
$F(I)\in\mathfrak{n}_{G}(B)$ or $F(I)\in \mathscr{I}(B)_{\mathscr{P}}$.
Suppose $F(I)\in \mathfrak{n}_{G}(B)$.
Then $F(I)$ is a node and there is an almost split sequence $0\ra F(I)\ra Q\ra {\rm Tr\,}D(F(I))\ra 0$ with $_BQ$ projective.
Since $F(\pi)$ is irreducible and $_BF(I)$ is indecomposable, by \cite[Chapter V, Theorem 5.3]{Aus1997}, we get $F(I/S)\in \add(_BQ)$.
Thus $F(I/S)$ is a projective $B$-module and $F(\pi)= 0$ in $\stmodcat{B}$. This is a contradiction and shows $F(I)\in \mathscr{I}(B)_{\mathscr{P}}$.

Similarly, we show that $G(J)$ lies in $\mathscr{I}(A)_{\mathscr{P}}$ for $J\in\mathscr{I}(B)_{\mathscr{P}}$.
By Lemma \ref{one-to-one}, $F:\mathscr{I}(A)_{\mathscr{P}}\ra \mathscr{I}(B)_{\mathscr{P}}$ and $F:\mathfrak{n}_{F}(A)\ra \mathfrak{n}_{G}(B)$ are bijections.

Let $P\in\mathscr{P}(A)_{\mathscr{I}}$ with $S$ as its top. By Lemma \ref{pi=stp}(3), $S$ is not projective. Now, Lemma \ref{exact-seq} implies $F'(P)\in \pmodcat{B}$. Thus  $F'(P)\in\bigtriangleup_A$ by Lemma \ref{one-to-one}, and therefore $F'(P)\in\mathscr{P}(B)_{\mathscr{I}}$. Similarly, $G'(Q)\in\mathscr{P}(A)_{\mathscr{I}}$ for $Q\in\mathscr{P}(B)_{\mathscr{I}}$. Thus Lemma \ref{one-to-one} yields the bijections $F':\mathscr{P}(A)_{\mathscr{I}}\ra \mathscr{P}(B)_{\mathscr{I}}$ and $F':\mathfrak{n}_{F}(A)\ra \mathfrak{n}_{G}(B)$. $\square$

\medskip
Suppose that $X$ is an $A$-module such that $A\stp=\add(_AX)$, where $A\stp$ stands for the full subcategory of $A\modcat$ consisting of all
$\nu$-stably projective $A$-modules.
By $\mbox{pre}(X)$ we denote the full subcategory of $A\modcat$ consisting of all those $A$-modules $M$ that have a minimal projective presentation
$P_1\ra P_0\ra M\ra 0$ with $P_1, P_0\in \add(X)$. Let $\Lambda:=\End_A(X).$ It follows from \cite[Chapter II, Proposition 2.5]{Aus1997} that the functor $\Hom_A(X,-):\mbox{pre}(X)\ra \Lambda\modcat$ is an equivalence of additive categories with a quasi-inverse $X\otimes_{\Lambda}-:\Lambda\modcat\ra \mbox{pre}(X)$. Let $R_X(M,N)$ be the $k$-submodule of $\Hom_A(M,N)$ consisting of all those homomorphisms of $A$-modules that factorize through a module in $\add(X)$. Then $R_X$ is an ideal of the $k$-category pre$(X)$. Note that $R_X(N,M)=R_A(N,M)$ for all $N\in A\modcat$ and $M\in \mbox{pre}(X)$ because a homomorphism $f:N\ra M$ of $A$-modules factorizes through a projective $A$-module must factorize through the projective cover $P_0\ra M$. We denote by $\underline{\mbox{pre}}(X)$  the quotient category of $\mbox{pre}(X)$ modulo the ideal $R_X$. Thus $\underline{\mbox{pre}}(X)$ is a full additive subcategory of $A\stmodcat$.

\medskip
{\bf Proof of Theorem \ref{main-thm}.}
Suppose that $F:\stmodcat{A}\ra \stmodcat{B}$ defines a stable equivalence between Artin algebras $A$ and $B$, where both algebras have positive $\nu$-dominant dimensions. Then $A\stp=\PI{A}$ and $B\stp=\PI{B}$ by Lemma \ref{pi=stp}(1). Let $_AX$ and $_BY$ be modules such that $A\stp=\add(_AX)$ and $B\stp=\add(_BY)$, respectively. Let $\Lambda:=\End(X)$ and $\Gamma:=\End_B(Y)$. Then $\Lambda$ and $\Gamma$ are the Frobenius parts of $A$ and $B$, respectively. In particular, $\Lambda$ and $\Gamma$ are self-injective algebras.

To show that $\Lambda$ and $\Gamma$ are stably equivalent, it is enough to show that $F$ induces an equivalence from $\underline{\mbox{pre}}(X)$ to $\underline{\mbox{pre}}(Y)$. Since $F$ is an equivalence, we need only to show that $F(M)$ lies in $\mbox{pre}(Y)$ for all $M\in \mbox{pre}(X)$.

In fact, take $M\in \mbox{pre}(X)$ and a minimal projective presentation: $P_1\ra P_0\ra M\ra 0$  with $P_1,P_0\in \add(_AX)$. We may assume that $M$ has no nonzero projective direct summands. Then the exact sequence
$$0\lra \Omega_A(M)\lra P_0\lra M\lra 0$$
is minimal and $\Omega_A(M)$ does not have any injective direct summands, that is, $\Omega_A(M)\in A\modcat_{\mathscr{I}}$. Since $P_1$ is a projective-injective $A$-module, we have $\Omega_A(M)\in A\modcat_{\mathscr{P}}$. So we write $\Omega_A(M)\simeq K_1\oplus K_2$ with $K_1\in \add(\bigtriangleup_A^c)$ and $K_2 \in \add(\mathfrak{n}_{F}(A))$. By Corollary \ref{exact}, we have a minimal exact sequence
$$(*)\quad 0\lra F(K_1)\oplus F'(K_2)\lra Q_0\lra F(M)\lra 0$$
in $B\modcat$
with $Q_0\in \PI{B}=\add(Y)$.

Next, we investigate the projective cover of $F(K_1)$. Let $P'_1$ be the projective cover of $K_1$. Then $P'_1\in \add(_AP_1)$ and $P'_1\in \PI{A}$. Note that $\Omega_A(K_1)$ lies in $A\modcat_{\mathscr{I}}$ and we can write $\Omega_A(K_1)= L_1\oplus L_2$ with $L_1\in \add(\bigtriangleup_A^c)$ and $L_2 \in \add(\Delta_A)$.  Applying Corollary \ref{exact} to the minimal exact sequence
$$0\lra \Omega_A(K_1)\lra P'_1\lra K_1\lra 0$$
in $A\modcat$, we get a minimal exact sequence of $B$-modules
$$0\lra F(L_1)\oplus F'(L_2)\lra Q'_1\lra F(K_1)\lra 0$$
with $Q'_1\in \PI{B}$.

Now, we investigate the projective cover of $F'(K_2)$. By Lemma \ref{bijective}, it follows from $K_2 \in \add(\mathfrak{n}_{F}(A))$ that $F'(K_2)\in \add(\mathfrak{n}_{G}(B))$, where $G$ is the quasi-inverse of the functor $F$. From the sequence $(*)$ and $Q_0\in \PI{B}$, we infer that the injective envelope of $F'(K_2)$ is projective. Since nodes are simple modules, it follows from Lemma \ref{pi=stp}(2) that the projective cover $_BQ_2'$ of $F'(K_2)$ is projective-injective. Thus  the minimal projective presentation of $F(M)$ is as follows.
$$ Q_1'\oplus Q_2'\lra Q_0\lra F(M)\lra 0$$
with $Q_1',Q_2', Q_0\in \PI{B}$. This yields $F(M)\in \mbox{pre}(Y)$. Similarly, we prove that the quasi-inverse $G$ of $F$ sends $N\in \mbox{pre}(Y)$ to $G(N)\in\mbox{pre}(X)$.

Finally, we reach the commutative diagram of functors in stable module categories:
$$\xymatrix{
\underline{\mbox{pre}}(X)\ar@<0.5ex>[r]^-{F.}\ar@{^{(}->}[d]
&\underline{\mbox{pre}}(Y)\ar@<0.5ex>[l]^-{G.}\ar@{^{(}->}[d] \\ \stmodcat{A}\ar@<0.5ex>[r]^-{F}& \stmodcat{B}\ar@<0.5ex>[l]^-{G}
}$$
where $F.$ and $G.$ stand for the restrictions of $F$ and $G$ to $\underline{\mbox{pre}}(X)$ and $\underline{\mbox{pre}}(Y)$ , respectively. It follows from the equivalence of $F$ that $F.$ is an equivalence of $k$-categories. $\square$

\medskip
The following example shows that the assumption of $\nu$-dominant dimensions in Theorem \ref{main-thm} cannot be dropped.
\begin{Bsp}{\rm
Let $A$ and $B$ be algebras given by the quivers with relations:
$$\begin{array}{ccc}
A:\xymatrix{
1\bullet\ar[r]^{\alpha} & \bullet 2\ar[d]^{\beta}\\
4\bullet\ar[u]_{\delta}&\bullet 3\ar[l]_{\gamma}}\\
\quad \beta\gamma=\delta\alpha=0,
\end{array}\qquad
\begin{array}{ccc}
B:\xymatrix{
1'\bullet\ar@<0.5ex>[r]^-{\alpha'}
&\bullet 2'\ar@<0.5ex>[l]^-{\beta'}
\quad
3'\bullet \ar@<0.5ex>[r]^-{\gamma'} & \bullet 4'\ar@<0.5ex>[l]^-{\delta'}}\\
\quad \beta'\alpha'=\delta'\gamma'=0.\\
\vspace{0.85cm}
\end{array}$$
We denote by $P(i)$ and $I(i)$ the indecomposable projective and injective modules corresponding to the vertex $i$, respectively. The indecomposable projective $A$-modules and $B$-modules are displayed, respectively.
$$
\xymatrix@R=.23cm@C=.01cm{
P(1)\\
1\ar@{-}[d]\\
2\ar@{-}[d]\\
3}\quad
\xymatrix@R=.23cm@C=.01cm{
P(2)\\
2\ar@{-}[d]\\
3}\quad
\xymatrix@R=.23cm@C=.01cm{
P(3)\\
3\ar@{-}[d]\\
4\ar@{-}[d]\\
1}\quad
\xymatrix@R=.23cm@C=.01cm{
P(4)\\
4\ar@{-}[d]\\
1}
\qquad\qquad
\xymatrix@R=.23cm@C=.01cm{
P(1')\\
1'\ar@{-}[d]\\
2'\ar@{-}[d]\\
1'}\quad
\xymatrix@R=.23cm@C=.01cm{
P(2')\\
2'\ar@{-}[d]\\
1'}\quad
\xymatrix@R=.23cm@C=.01cm{
P(3')\\
3'\ar@{-}[d]\\
4'\ar@{-}[d]\\
3'}\quad
\xymatrix@R=.23cm@C=.01cm{
P(4')\\
4'\ar@{-}[d]\\
3'}
$$
The indecomposable injective $A$-modules and $B$-modules are given as follows.
$$
\xymatrix@R=.23cm@C=.01cm{
I(1)\\
3\ar@{-}[d]\\
4\ar@{-}[d]\\
1}\quad
\xymatrix@R=.23cm@C=.01cm{
I(2)\\
1\ar@{-}[d]\\
2}\quad
\xymatrix@R=.23cm@C=.01cm{
I(3)\\
1\ar@{-}[d]\\
2\ar@{-}[d]\\
3}\quad
\xymatrix@R=.23cm@C=.01cm{
I(4)\\
3\ar@{-}[d]\\
4}
\qquad\qquad
\xymatrix@R=.23cm@C=.01cm{
I(1')\\
1'\ar@{-}[d]\\
2'\ar@{-}[d]\\
1'}\quad
\xymatrix@R=.23cm@C=.01cm{
I(2')\\
1'\ar@{-}[d]\\
2'}\quad
\xymatrix@R=.23cm@C=.01cm{
I(3')\\
3'\ar@{-}[d]\\
4'\ar@{-}[d]\\
3'}\quad
\xymatrix@R=.23cm@C=.01cm{
I(4')\\
3'\ar@{-}[d]\\
4'}
$$
Then $A\stp=\add\big(P(1)\oplus P(3)\big)$, $B\stp=\add\big(P(1')\oplus P(3')\big)$ and
$\ndd(A)=\ndd(B)=2$. The Frobenius parts $\Lambda$ and $\Gamma$ of $A$ and $B$ are given by the quivers with relations, respectively.
$$\begin{array}{ccc}
\Lambda:\quad \xymatrix{
1\bullet \ar@<0.5ex>[r]^{\alpha} &\bullet 3\ar@<0.5ex>[l]^{\gamma}}
&\hspace{2cm}&
\Gamma:\quad\xymatrix{
1'\bullet\ar@(ru,rd)^{\alpha'}}
\quad
\xymatrix{
3'\bullet\ar@(ru,rd)^{\gamma'}}\\
\qquad \quad \alpha\gamma=\gamma\alpha=0,
&\hspace{2cm}&
\qquad\alpha'^2= \gamma'^2=0.
\end{array}$$
It follows from\cite[Theorem 2.10]{MV1980} (see Lemma \ref{node}(2) below) that both $A$ and $B$ are stably equivalent to the path algebra $C$ of the quiver
$$C:\xymatrix@C=0.6cm{1\bullet \ar[r]& \bullet\ar[r]^(0){2}& \bullet 5}
\quad
\xymatrix@C=0.6cm{3\bullet \ar[r]&\bullet\ar[r]^(0){4}\ar[r]&\bullet 6}.$$
Thus $A$ and $B$ are stably equivalent, and so are $\Lambda$ and $\Gamma$ by Theorem \ref{main-thm}. Now, we consider the stably equivalent algebras $A$ and $C$. Clearly, $\ndd(C)=0$ and the Frobenius part of $C$ is $0$. Thus the Frobenius part $\Lambda$ of $A$ is not stably equivalent to the Frobenius part of $C$. This shows that the assumption on $\nu$-dominant dimensions of Artin algebras in Theorem \ref{main-thm} cannot be omitted. Observe that $\dd(C)=1$. This shows that the $\nu$-dominant dimensions in Theorem \ref{main-thm}  cannot be weakened to dominant dimensions either.
}\end{Bsp}

\section{Auslander--Reiten conjecture for two classes of algebras\label{sect5}}

This section is devoted to the proof of Theorem \ref{cma}, namely we show that (ARC) holds true for principal centralizer matrix algebras over arbitrary fields and Frobenius-finite algebras over algebraically closed fields. Though the two large classes of algebras are rather different, they share a common feature that Frobenius parts are representation-finite.

\subsection{Auslander--Reiten conjecture for Frobenius-finite algebras}

Recall that the Auslander--Reiten conjecture states that stably equivalent algebras have the same number of non-projective, non-isomorphic simple modules. The following result establishes  a relation of validity of the conjecture between algebras and their Frobenius parts.

\begin{Lem}\label{alg to F-part}
Let $A$ and $B$ be stably equivalent Artin algebras, and let $\Lambda$ and $\Gamma$ be the Frobenius parts of $A$ and $B$, respectively.

$(1)$ If $A$ and $B$ have no semisimple direct summands, then neither  do $\Lambda$ and $\Gamma$.

$(2)$ Suppose that $A$ and $B$ have no nodes. Then $\Lambda$ and $\Gamma$ are stably equivalent. If, in addition, $\Lambda$ and $\Gamma$ have the same number of non-isomorphic, non-projective simples, then so do $A$ and $B$.

$(3)$ Suppose that $\ndd(A)\ge 1$ and $\ndd(B)\ge 1$. If one of $\Lambda$ and $\Gamma$ is a Nakayama algebra, then $A$ and $B$ have the same number of non-isomorphic, non-projective simples.
\end{Lem}

{\it Proof.} Let $X\in A\modcat$ and $Y\in B\modcat$ such that $A\stp=\add(_AX)$  and $B\stp=\add(_BY)$, and let $\Lambda:=\End(_AX)$ and $\Gamma:=\End_B(Y)$. Then  $\Lambda$ and $\Gamma$ are the Frobenius parts of $A$ and $B$, respectively, and therefore they are self-injective Artin algebras.

(1) We show that if $\Lambda$ has semisimple direct summands then so does $A$. Indeed, without loss of generality, we may assume that $A$ is a basic algebra and $_AX$ is a basic $A$-module. Then $\Lambda$ is a basic algebra. Since $\Lambda$ has semisimple direct summands, there is a nonzero central idempotent $e$ of $\Lambda$ such that $e\Lambda e$ is semisimple and $\Lambda=e\Lambda e\times (1-e)\Lambda (1-e)$. In particular, $(1-e)\Lambda e=e\Lambda (1-e)=0$. Note that $e\Lambda e$ is basic. It follows from the Wedderburn--Artin theorem that $e\Lambda e$ is isomorphic to a product of finitely many division rings. Let $e_0$ be a primitive idempotent of $\Lambda$ with $e_0\in e\Lambda e$. Then $e_0\Lambda e_0$ is a division ring and $\Lambda= e_0\Lambda e_0\times (e-e_0)\Lambda (e-e_0)\times (1-e)\Lambda (1-e)$. Particularly, $(1-e_0)\Lambda e_0=e_0\Lambda (1-e_0)=0$. Since the evaluation functor $\Hom_A(X,-):A\modcat\ra \Lambda\modcat$ induces an equivalence $\add(_AX)\simeq \Lambda\pmodcat$ of additive categories, there is an indecomposable summand $X_0$ of $_AX$ such that $\Hom_{A}(X,X_0)\simeq \Lambda e_0$. Then $\End_{A}(X_0)\simeq e_0\Lambda e_0$ is a division ring, $\Hom_{A}(X/X_0,X_0)=(1-e_0)\Lambda e_0=0$, and $\Hom_{A}(X_0,X/X_0)=e_0\Lambda (1-e_0)=0$. As $_AX\in A\stp$ is basic, $\top(_AX)\simeq \soc(_AX)$. Thus $\top(_AX_0)\simeq \soc(_AX_0)$. Since $\End_A(X_0)$ is a division ring, $X_0$ must be a simple $A$-module in $\PI{A}$. Then $\Hom_A(X_0,P)=\Hom_A(P,X_0)=0$ for any indecomposable projective $A$-module $P$ which is not isomorphic to $X_0$. Thus $A\simeq \End_A(X_0)\times \End_A(A/X_0)$. In particular, $\End_A(X_0)$ is a semisimple direct summand of $A$.

(2) Since we are only concerned with non-projective simple modules, we may assume that $A$ and $B$ have no semisimple direct summands. It follows from \cite[Theorem 2.6]{MV1990} that $\Lambda$ and $\Gamma$ are stably equivalent.

Assume further that $\Lambda$ and $\Gamma$  have the same number of non-isomorphic, non-projective simple modules. We show that $A$ and $B$ have the same number of non-isomorphic, non-projective simple modules. Indeed, it follows from \cite[Lemma 2.5]{MV1990} which holds true also for Artin algebras, that $A$ and $B$ have the same number of non-isomorphic, non-projective, simple modules whose projective covers are not $\nu$-stably projective. It remains to show that $A$ and $B$ have the same number of non-isomorphic, non-projective, simple modules whose projective covers are $\nu$-stably projective. Note that a projective simple module is not $\nu$-stably projective. Otherwise, it would be a projective-injective simple module, and therefore $A$ and $B$ would have semisimple direct summands. Thus we have to show that $A$ and $B$ have the same number of non-isomorphic simple modules whose projective covers are $\nu$-stably projective. As $A\stp=\add(_AX)$ and $B\stp=\add(_BY)$, we need to show that $\Lambda$ and $\Gamma$ have the same number of non-isomorphic simple modules. Note that $\Lambda$ and $\Gamma$ do not have projective simple modules by (1). By assumption, $\Lambda$ and $\Gamma$ have the same number of non-isomorphic simple modules. Hence $A$ and $B$ have the same number of non-isomorphic, non-projective simple modules.

(3) Without loss of generality, we assume that $A$ and $B$ have no semisimple direct summands. We have to show that $A$ and $B$ have the same number of non-isomorphic, non-projective simples. Indeed, due to $\ndd(A)\ge 1$ and $\ndd(B)\ge 1$, it follows from Lemma \ref{bijective} that $A$ and $B$ have the same number of non-isomorphic, non-projective, simple modules whose projective covers are not injective. By Lemma \ref{pi=stp}(1), $A\stp=\PI{A}$ and $B\stp=\PI{B}$. It remains to show that $A$ and $B$ have the same number of non-isomorphic, non-projective, simple modules whose projective covers are $\nu$-stably projective. By Theorem \ref{main-thm}, $\Lambda$ and $\Gamma$ are stably equivalent.
Assume that one of $\Lambda$ and $\Gamma$ is a Nakayama algebra.
By \cite[Theorem 1.3]{Reiten78} which says that if an Artin algebra is stably equivalent to a Nakayama algebra then the two algebras have the same number of non-isomorphic, non-projective simple modules, we deduce that $\Lambda$ and $\Gamma$ have the same number of non-isomorphic, non-projective simples. An argument similar to the proof of (2) shows that $A$ and $B$ have the same number of non-isomorphic, non-projective, simple modules whose projective covers are $\nu$-stably projective. Thus $A$ and $B$ have the same number of non-isomorphic, non-projective simples.
$\square$

\medskip
A finite-dimensional $k$-algebra $A$ over a field $k$ is called a {\em Morita algebra} if $A$ is isomorphic to \mbox{$\End_H(H\oplus M)$} for $H$ a finite-dimensional self-injective $k$-algebra and $M$ a finitely generated $H$-module \cite{KY13}. If $H$ is symmetric, then the Morita algebra $A$ is called a {\em gendo-symmetric algebra} \cite{FK16}. In this case, the Frobenius part of $A$ is Morita equivalent to $H$. Recently, it has been shown that $S_n(c,k)$ is always a gendo-symmetric algebra \cite[Theorem 1.1(2)]{Xi2022}. An algebra $A$ is a Morita algebra if and only if $\ndd(A)\ge 2$ by \cite[Proposition 2.9]{fhk21}.

For $c\in M_n(k)$, we denote by $k[c]$ the unitary subalgebra of $M_n(k)$ generated by $c$. Let $\varphi:k[x]\to k[c]$ be the surjective homomorphism of algebras, defined by $x\mapsto c$. Then $\Ker(\varphi)=(m_c(x))$ where $m_c(x)$ is the minimal polynomial of $c$ over $k$, and $\varphi$ induces an isomorphism $\bar{\varphi}: k[x]/(m_c(x))\simeq k[c]$ of algebras. Let $A_c:=k[x]/(m_c(x))$, and let $k^n$ be the $n$-dimensional vector space over $k$ consisting of column vectors. Then $k^n$ is naturally a $k[c]$-module, and therefore an $A_c$-module via $\bar{\varphi}$. By definition, $S_n(c,k)^{\opp}\simeq \End_{A_c}(k^n)$. If we write $m_c(x):=\prod^{s}_{i=1} f_i(x)^{n_i}$ with all $f_i(x)$ pairwise coprime irreducible polynomials and set $B_i:=k[x]/(f_i(x)^{n_i})$ for $1\le i\le s$, then it follows from the Chinese remainder theorem that $A_c:=k[x]/(m_c(x))\simeq  \prod_{i=1}^{s}B_i$. Now, we decompose the $A_c$-module $k^n=\bigoplus^{s}_{i=1} M_i$ such that $M_i$ is a direct sum of representatives of the isomorphism classes of indecomposable direct summands of $k^n$ lying in the block $B_i$. Then $S_n(c,k)^{\opp}\simeq \prod_{i=1}^{s}\End_{B_i}(M_i).$
 Clearly, $k^n$ is a faithful $M_n(k)$-module and $k[c]$ is a subalgebra of $M_n(k)$. Thus $k^n$ is also a faithful $k[c]$-module. This implies that $M_i$ is a faithful $B_i$-module for $1\le i\le s$. As $B_i$ is a symmetric Nakayama algebra (see \cite[Section V.1 Example, pp. 140-141]{Aus1997}), we know that $M_i$ is a generator for $B_i\modcat$ and $\End_{B_i}(M_i)$ is a gendo-symmetric algebra for $1\le i\le s$. Due to the isomorphisms $S_n(c,k)\simeq S_n(c',k)\simeq S_n(c,k)^{\opp}$ as algebras, where $c'$ is the transpose of the matrix $c$, we see that the gendo-symmetric algebra $S_n(c,k)$ has its Frobenius part Morita equivalent to $B_i$ for $1\le i\le s$. Thus $S_n(c,k)$ is a gendo-symmetric algebra such that its Frobenius part is a symmetric Nakayama algebra.

\medskip
{\bf Proof of Theorem \ref{cma}}. (1) Let $c\in M_n(k)$ and $d\in M_m(k)$. Suppose that $S_n(c,k)$ and $S_m(d,k)$ are stably equivalent. Since $S_n(c,k)$ and $S_m(d,k)$ are gendo-symmetric, it follows from \cite[Proposition 2.9]{fhk21} that $\ndd(S_n(c,k))\ge 2$ and $\ndd(S_m(d,k))\ge 2$. Thanks to Theorem \ref{main-thm}, the Frobenius parts of both $S_n(c,k)$ and $S_m(d,k)$ are also stably equivalent. Note that the Frobenius parts of both $S_n(c,k)$ and $S_m(d,k)$ are Nakayama algebras. It follows from Lemma \ref{alg to F-part}(3) that $S_n(c,k)$ and $S_m(d,k)$ have the same number of non-isomorphic, non-projective simples.

(2) Assume that $A$ and $B$ are Artin $k$-algebras over a commutative Artin ring $k$. Given a stable equivalence between $A$ and $B$, we get a stable equivalence between $A'$ and $B'$ both of which have no nodes.  Let $\Lambda'$ and $\Gamma'$ be the Frobenius parts of $A'$ and $B'$, respectively. Then $\Lambda'$ and $\Gamma'$ are stably equivalent by Lemma \ref{alg to F-part}(2).

Now, assume that $k$ is an algebraically closed field and that $A$ is Frobenius-finite. Then $A'$ is Frobenius-finite by Lemma \ref{node}(4), that is, $\Lambda'$ is representation-finite and therefore $\Gamma'$ is representation-finite. Since Auslander--Reiten conjecture holds true for a stable equivalence between representation-finite $k$-algebras over an algebraically closed field $k$ (see \cite[Theorem 3.4]{MV1985}), $\Lambda'$ and $\Gamma'$ have the same number of non-isomorphic, non-projective simple modules. By Lemma \ref{alg to F-part}(2), $A'$ and $B'$ have the same number of non-isomorphic non-projective simple modules, and therefore $A$ and $B$ have the same number of non-isomorphic non-projective simple modules by Lemma \ref{node}(3).
$\square$

\medskip
The following result is an immediate consequence of Theorems \ref{main-thm} and \ref{cma}(2). Here algebras considered may have nodes.

\begin{Koro}\label{main-cor}
Every stable equivalence of Morita $k$-algebras over a field $k$ induces a stable equivalence of their Frobenius parts. In particular, if $A$ and $B$ are stably equivalent Morita algebras over an algebraically closed field and if one of $A$ and $B$ is Frobenius-finite, then the Auslander--Reiten conjecture holds true for $A$ and $B$, namely they have the same number of non-isomorphic, non-projective simples.
\end{Koro}

\subsection{A conjecture}
As we know, derived equivalences do not have to preserve the delooping levels of algebras in general (see Remark \ref{rmk-tilt}).
This happens also for global and finitistic dimensions of algebras. However, it is well known that derived equivalences preserve finiteness of global and finitistic dimensions. Also, derived equivalences preserve finiteness of $\phi$- and $\psi$-dimensions of algebras (see \cite{flm}). All of these phenomena suggest the following conjecture.

\smallskip
{\bf Conjecture.}  If $A$ and $B$ are derived equivalent noetherian rings, then $\del(A)<\infty$ if and only if $\del(B)<\infty,$ that is, the finiteness of delooping levels of algebras is invariant under derived equivalences.

\medskip
The validity of this conjecture can be applied to adjudge finiteness of finitistic dimensions of derived equivalent algebras if the opposite algebra of one of the algebras has finite delooping level. It can also be used to test whether two concrete algebras are derived equivalent or not by calculating their delooping levels.

\medskip
{\footnotesize

Changchang Xi, School of Mathematical Sciences, Capital Normal University, 100048 Beijing, P. R. China
and School of Mathematics and Statistics, Shaanxi Normal University, 710119 Xi'an, P. R. China

{\tt Email: xicc@cnu.edu.cn (C.C.Xi)}

\medskip
Jinbi Zhang, School of Mathematical Sciences, Peking University, 100871 Beijing, P. R. China

Current address:
School of Mathematical Sciences, Anhui University,  230601 Hefei, P. R. China

{\tt Email: zhangjb@ahu.edu.cn (J.B.Zhang)}
}

\end{document}